\documentclass[11pt, dina4]{article}
\usepackage{amssymb, amsmath,amsthm}
\usepackage{graphicx}
\usepackage{multirow}
\usepackage{epstopdf}
\usepackage{float}
\usepackage{diagbox}
\usepackage{subfig}
\usepackage{soul}
\usepackage{color}
\usepackage[margin=1in]{geometry}

\newcommand{\p}{\partial}
\newcommand{\og}{\omega}
\newcommand{\Og}{\Omega}
\newcommand{\fl}[2]{\frac{#1}{#2}}
\newcommand{\dt}{\delta}

\newcommand{\nn}{\nonumber}
\newcommand{\ap}{\alpha}
\newcommand{\bt}{\beta}

\newcommand{\veps}{\varepsilon}

\newcommand{\Dt}{\Delta}

\newcommand{\be}{\begin{equation}}
\newcommand{\ee}{\end{equation}}
\newcommand{\ba}{\begin{array}}
\newcommand{\ea}{\end{array}}
\newcommand{\bea}{\begin{eqnarray}}
\newcommand{\eea}{\end{eqnarray}}
\newcommand{\beas}{\begin{eqnarray*}}
\newcommand{\eeas}{\end{eqnarray*}}
\newtheorem{remark}{Remark}[section]
\newtheorem{lemma}{Lemma}[section]

\newcommand{\bx}{{\bf x} }
\newcommand{\by}{{\bf y} }

\newcommand{\cc}{\color{blue}}
\newcommand{\bb}{\vskip 10pt}
\definecolor{ForestGreen}{rgb}{0.0, 0.5, 0.0}

\title{A  universal solution scheme for fractional and classical PDEs}
\author{Yixuan Wu\thanks{Department of Mathematics and Statistics, Missouri University of Science and Technology, Rolla, MO 65409 (Email:  ywx7c@mst.edu)}, \ \
Yanzhi Zhang\thanks{Department of Mathematics and Statistics, Missouri University of Science and Technology, Rolla, MO 65409 (Email:  zhangyanz@mst.edu; URL:  {http://web.mst.edu/$\sim$zhangyanz})}}
\begin{document}
\date{}
\maketitle

\begin{abstract} 
We propose a unified meshless method   to solve classical and fractional PDE problems with $(-\Dt)^{\fl{\ap}{2}}$ for $\ap \in (0, 2]$. 
The classical ($\ap = 2$) and fractional ($\ap < 2$) Laplacians, one local and the other nonlocal, have distinct properties. 
Therefore, their numerical methods and computer implementations are usually incompatible. 
We notice that for any $\ap \ge 0$, the Laplacian $(-\Dt)^{\fl{\ap}{2}}$ of generalized inverse multiquadric  (GIMQ) functions  can be analytically written by the Gauss hypergeometric function, and thus propose a GIMQ-based method. 
Our method unifies the discretization of classical and fractional Laplacians and also bypasses numerical approximation to the hypersingular integral of fractional Laplacian. 
These two merits distinguish our method from other existing methods for the fractional Laplacian. 
 Extensive numerical experiments are carried out to test the performance of our method.
 Compared to other methods, our method can achieve high accuracy with fewer number of unknowns, which effectively reduces the storage and computational requirements in simulations of fractional PDEs. 
 Moreover, the meshfree nature makes it free of geometric constraints and enables simple implementation for any dimension $d \ge 1$. 
Additionally, two approaches of selecting shape parameters, including condition number-indicated method and random-perturbed  method, are studied to avoid the ill-conditioning issues when large number of points. 
\end{abstract}


\section{Introduction}
\setcounter{equation}{0}
\label{section1}

Over the last decade, fractional partial differential equations (PDEs) have been well recognized for their ability to describe anomalous diffusion phenomena in many complex systems \cite{Cusimano2015, delCastillo2003, Hanert2011}. 
Mathematically, anomalous diffusion can be modeled by the fractional Laplacian $(-\Dt)^{\fl{\ap}{2}}$ (for $\ap < 2$), in contrast to the classical Laplacian $\Dt = \p_{xx} + \p_{yy} + \p_{zz}$ for (normal) diffusion. 
It is well-known that the fractional Laplacian $(-\Dt)^{\fl{\ap}{2}}$ collapses to the classical Laplacian $-\Dt$ as $\ap \to 2^{-}$, but the properties of these two operators are essentially different. 
The fractional Laplacian $(-\Dt)^\fl{\ap}{2}$, representing the infinitesimal generator of a symmetric $\ap$-stable L\'evy process, is a nonlocal operator, while the classical Laplacian  $\Dt$  is a local operator. 
Due to this distinct difference, most existing numerical methods developed for the classical Laplacian can not be applied to solve problems with the fractional Laplacian. 
So far,  numerical methods for the fractional Laplacian  $(-\Dt)^{\fl{\ap}{2}}$ still remain limited. 
Moreover, some discretization techniques (e.g., finite element methods) for the classical and fractional Laplacians are incompatible, and different computer codes are required for their practical implementation. 
In this work, we propose a meshless pseudospectral method based on the generalized inverse multiquadric functions to solve normal and anomalous diffusion problems. 
It not only enriches the collection of numerical methods for studying the fractional Laplacian, but also provides a unified approach to solve both classical and fractional PDEs.

Let $\Og \subset {\mathbb R}^d$ (for $d = 1, 2$, or $3$) be an open bounded domain.  
Consider the following diffusion problem  with  Dirichlet boundary conditions: 
\bea
\label{diffusion}
\p_t u(\bx, t) = -\kappa(-\Dt)^{\fl{\ap}{2}}u + f(\bx, t, u), &\, &\mbox{for} \ \, \bx \in \Og, \ t > 0, \\
\label{diffusion-BC}
u(\bx, t) = g(\bx, t), && \mbox{for} \ \, \bx \in \Upsilon, \  t \ge 0, 
\eea
where the diffusive power $\ap \in (0, 2]$, and diffusion coefficient $\kappa > 0$.  
The notation $\Upsilon$ depends on $\ap$, i.e., $\Upsilon = \p{\Og}$ for $\ap = 2$,\, while $\Upsilon = \Og^c = {\mathbb R}^d\backslash\Og$ if $\ap < 2$. 
The Laplace operator $(-\Dt)^{\fl{\ap}{2}}$ is defined via an $\ap$-parametric pseudo-differential form \cite{Landkof,Samko}:
\begin{eqnarray}
\label{pseudo}
(-\Delta)^{\fl{\alpha}{2}}u({\bx}) = \mathcal{F}^{-1}\big[|\xi|^\alpha \mathcal{F}[u]\big], \qquad \mbox{for} \ \  \ap > 0, 
\end{eqnarray}
where $\mathcal{F}$ represents the Fourier transform with the associated inverse transform $\mathcal{F}^{-1}$. 
The pseudo-differential operator in \eqref{pseudo} gives a unified definition to the classical and fractional Laplacians. 
For $\ap = 2$,  it reduces to the spectral representation of the classical Laplacian $-\Dt$.  
Thus, (\ref{diffusion})--(\ref{diffusion-BC}) collapses to a classical (normal) diffusion problem with Dirichlet boundary conditions on $\p\Og$. 
While $\ap \in (0, 2)$, the operator \eqref{pseudo} defines the fractional Laplacian that models anomalous diffusion due to L\'evy flights, 
and in this case the Dirichlet boundary condition is extended to $\Og^c$ (i.e., the complement of domain $\Og$). 

The fractional Laplacian $(-\Delta)^{\fl{\alpha}{2}}$ can be also defined via a hypersingular integral  \cite{Landkof, Samko}, i.e.,
 \bea\label{integralFL}
(-\Delta)^{\fl{\alpha}{2}}u(\bx) = C_{d,\alpha}\,{\rm P. V.}\int_{{\mathbb R}^d}\fl{u(\bx) - u({\bf y})}{|\bx -{\bf y}|^{d+\ap}}d{\bf y}, \qquad \mbox{for} \ \  \ap \in (0, 2),
\eea
where ${\rm P. V.}$ stands for the principal value integral, and the normalization constant 
\beas
C_{d,\alpha}=\fl{2^{\ap-1} \ap\,\Gamma({\ap+d}/{2})}{\sqrt{\pi^{d}}\,\Gamma(1 -{\ap}/{2})}
\eeas
with $\Gamma(\cdot)$ being the Gamma function. 
This integral operator provides a pointwise definition of the fractional Laplacian. 
Note that the pseudo-differential operator in (\ref{pseudo}) unifies the definition of classical and fractional Laplacians so as to provide a foundation for developing compatible schemes for these two operators, but  it is challenging to incorporate non-periodic boundary conditions into (\ref{pseudo}). 
In contrast to (\ref{pseudo}), the pointwise definition of the fractional Laplacian in \eqref{integralFL} can easily incorporate  other boundary conditions, but it is incompatible to the classical Laplacian (i.e., $\ap \neq 2$ in \eqref{integralFL}). 
As shown in \cite{Samko, Kwasnicki2017},  the two definitions of the fractional Laplacian in \eqref{pseudo} and \eqref{integralFL} are equivalent for functions in the Schwartz space. 
Hence, we want to combine the advantages of both definitions and develop a numerical method that is compatible for both classical and fractional Laplacians and can easily work with Dirichlet boundary conditions as in (\ref{diffusion-BC}). 

So far, numerical methods for the fractional Laplacian $(-\Dt)^{\fl{\ap}{2}}$ still remain scant. 
Most of them are  based on the pointwise definition of the fractional Laplacian and thus need to approximate the hypersingular integral in \eqref{integralFL}; see \cite{Acosta2017,   Ainsworth2018B, Duo2018, Acosta2019,  Duo-FDM2019, Duo-TFL2019, Bonito2019} and references therein. 
On the other side, meshless RBF-based methods have been recognized for their flexiblity with complex geometries, high accuracy and efficiency,  and simplicity of implementation. 
They have been widely applied to solve classical PDEs \cite{Fornberg2015,Fornberg2015b,Bhatia2016}. 
But, the application of RBF-based methods in solving nonlocal or fractional PDEs is still very recent. 
To the best of our knowledge, two meshless RBF-based methods were recently proposed  in the literature  to solve problems governed by the fractional Laplacian $(-\Dt)^{\fl{\ap}{2}}$ (for $\ap < 2$) -- one is the Wendland RBF method in \cite{Rosenfeld2019}, and the other is the Gaussian RBF method in \cite{Burkardt0020}. 
The Wendland RBF-based method directly approximates the Fourier integrals in definition \eqref{pseudo} and considers the extended Dirichlet boundary conditions (for $\ap < 2$) on a large truncated region $\og \subset \Og^c$. 
It then solves the problem on $\Og\, \cup\, \og$ (instead of $\Og$), and consequently requires significant amount of RBF points and demands high computational cost. 
In contrast, the Gaussian RBF-based method integrates the extended boundary conditions into the scheme via the pointwise definition (\ref{integralFL}), and no boundary truncation is needed. 
Moreover, it can discretize both classical and fractional Laplacians in a single scheme \cite{Burkardt0020}.

The aim of this work is to develop a unified meshless pseudospectral method to solve both classical ($\ap = 2$) and fractional ($\ap < 2$) PDEs. 
Its unique feature -- compatibility with the classical Laplacian -- makes our method distinguish from other numerical methods (e.g., in \cite{Duo2018, Duo-FDM2019, Ainsworth2018B, Acosta2019,Bonito2019, Acosta2017,  Bond2015, Bonito2019}) for the fractional Laplacian. 
Moreover, our method takes great advantage of the Laplacian $(-\Dt)^{\fl{\ap}{2}}$ (for both $\ap = 2$ and $\ap < 2$) of generalized inverse multiquadric functions so as to bypass numerical approximations to the hypersingular integral of fractional Laplacian in (\ref{integralFL}).  
Extensive numerical experiments are carried out to test the performance of our method. 
Compared to other methods, our method can achieve high accuracy with fewer number of unknowns, which effectively reduces the storage and computational requirements in simulations of fractional PDEs. 
Moreover, the meshfree nature makes it free of geometric constraints and enables simple implementation  for any dimension $d \ge 1$. 
We also study two approaches, namely the condition number-indicated approach and random-perturbed approach, for selecting the shape parameters. 
Numerical studies show that they can effectively control the ill-conditioning issues and improve the performance of our method. 

The paper is organized as follows.  
In Section \ref{section2}, we introduce radial basis functions and the properties of generalized (inverse) multiquadric functions. 
In Section \ref{section3},  we propose a meshless method based on generalized inverse multiquadric functions 
 to solve the diffusion problem \eqref{diffusion}--\eqref{diffusion-BC}. 
In Section \ref{section4}, we test the accuracy and compare the proposed method with the Gaussian RBF-based method in \cite{Burkardt0020}. 
Two approaches in selecting the shape parameter are also studied. 
In Section \ref{section5},  we further test the performance of our method in solving elliptic problems and diffusion equations.  
Finally, discussion and conclusion are made in Section \ref{section6}.

\section{Radial basis functions}
\label{section2}
\setcounter{equation}{0}

Radial basis functions (RBFs) are a family of functions that depend on the distance of point $\bx$ to a given center point $\by$, i.e. $\varphi(\bx)=\varphi(|\bx- \by|)$ for $\bx,\, \by \in {\mathbb R}^d$. 
RBFs have been well recognized for their success in interpolating high-dimensional  scattered data. 
The application of RBFs in solving PDEs was first proposed by  Kansa in 1990 \cite{Kansa90I, Kansa90II}. 
Since then,  RBF-based methods have been widely applied to solve problems arising in various applications (see \cite{Fornberg2015,Fornberg2015b,Bhatia2016} and references therein). 

In the literature, RBFs are usually divided into two main categories:  globally-supported functions (e.g., see Table \ref{Tab0}) and compactly-supported functions (e.g., Wendland function). 
\begin{table}[htb!]
\begin{center}
 \begin{tabular}{|lcr|} 
 \hline
 Name &\qquad\qquad   & Definition $\varphi(r)$  \\ 
 \hline
 Gaussian     && ${e^{-(\varepsilon r)}}^2$ \\
 Multiquadric && $\sqrt{1 + (\varepsilon r)^2}$  \\
 Inverse multiquadric & & $\fl{1}{\sqrt{1 + (\varepsilon r)^2}}$\\
 Polyharmonic spline & & $r^{2m}{\rm ln}(r),\ m\in \mathbb{N}$  \\ 
 \hline
\end{tabular}
\caption{Some examples of globally-supported  RBFs with $r = |\bx - \by|$. }\label{Tab0}
\end{center}
\end{table}
The Gaussian and multiquadric functions are two well-known globally supported RBFs, and both of them are infinitely differentiable. 
The multiquadric function  was first proposed by Hardy to interpolate scattered data on a topographical map \cite{Hardy1971}. 
Later, Franke compared the multiquadric interpolation with other methods and concluded that  the multiquadric  function outperforms others in many aspects, including  numerical accuracy, computational time, and implementation simplicity \cite{Franke1979, Franke1982}.  
Recently, a broad class of generalized multiquadric functions have received a lot of attention in the literature \cite{Chenoweth2009, Sarra2009}.

\subsection{Generalized multiquadrics RBFs}
\label{section2-1}

The generalized multiquadric radial basis functions take the following form \cite{Chenoweth2009, Sarra2009}: 
\bea\label{fun2-1-1}
\varphi(r) = \big(1+\varepsilon^2r^2\big)^\beta, \qquad \mbox{for} \ \ \bt \in {\mathbb R}\backslash {\mathbb N}^0, 
\eea
where $\varepsilon > 0$ denotes the shape parameter, and ${\mathbb N}^0$ represents the set of non-negative integers. 
The function in  (\ref{fun2-1-1}) covers a wide range of infinitely differentiable RBFs, including, e.g.,
\begin{itemize}\itemsep -1pt
\item the well-known Hardy's multiquadric function when $\bt = \fl{1}{2}$ \cite{Hardy1971, Hardy1990};
\item the inverse multiquadric function when $\beta = -\fl{1}{2}$ \cite{Ku2020,Soleymani2018};
\item the inverse quadratic function when $\bt = -1$.
\end{itemize}
For convenience of discussion, we will refer function \eqref{fun2-1-1} with $\bt > 0$ as the {\it generalized multiquadric} (GMQ) function, while the {\it generalized inverse multiquadric} (GIMQ) function if $\bt < 0$. 
It shows that  the GIMQ functions are strictly positive definite, while  the GMQ functions are strictly conditionally positive definite of order $\lceil \bt \rceil$, where $\lceil \cdot \rceil$ denotes the ceiling function \cite{Chenoweth2009,Sarra2009}. 
The strictly positive definite is an important property to ensure the system matrix from interpolation to be invertible. 

\begin{figure}[htb!]
\centerline{
\includegraphics[height = 5.86cm, width = 8.46cm]{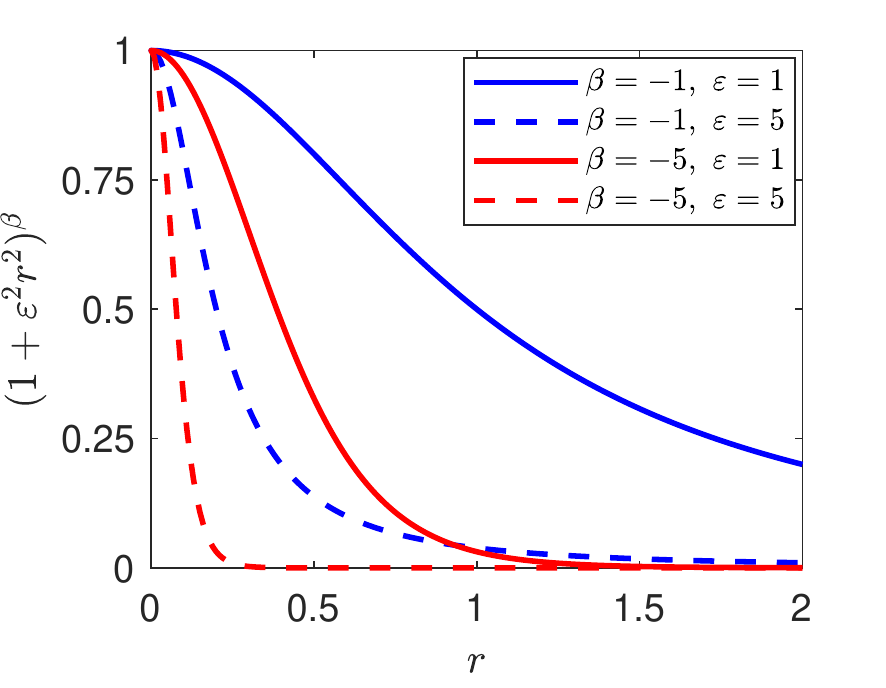} }
\caption{Illustration of GIMQ functions for different shape parameter $\veps$ and power $\bt$.}
\label{Fig2-1-1}
\end{figure}

Note that both GIMQ and Gaussian RBFs are strictly positive definite. 
The shape of Gaussian functions is controlled by  parameter $\veps$ --  the smaller the shape parameter $\veps$, the flatter the Gaussian function. 
While the GIMQ functions are controlled by both shape parameter $\veps$ and power $\bt$; see the comparison of their effects in Fig. \ref{Fig2-1-1}. 
It shows that the GIMQ functions monotonically decrease as $r \to \infty$. 
The larger the power $\bt$ (or the smaller the shape parameter $\veps$), the flatter the GIMQ function. 
The effect of $\veps$ is more significant in determining the shape of GIMQ functions.

\begin{figure}[htb!]
\centerline{\includegraphics[height = 5.76cm, width = 7.86cm]{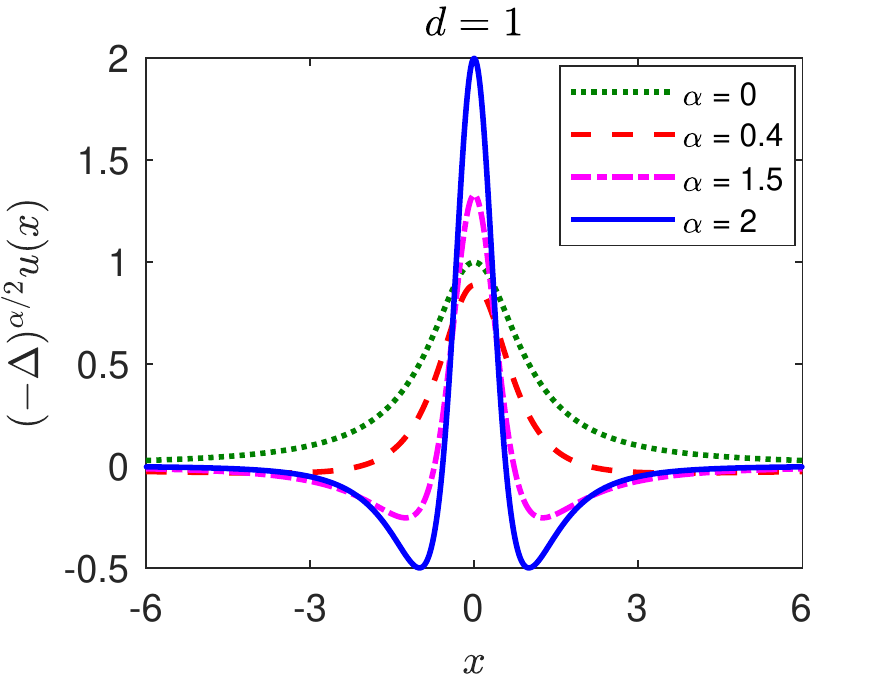} 
\includegraphics[height = 5.76cm, width = 7.86cm]{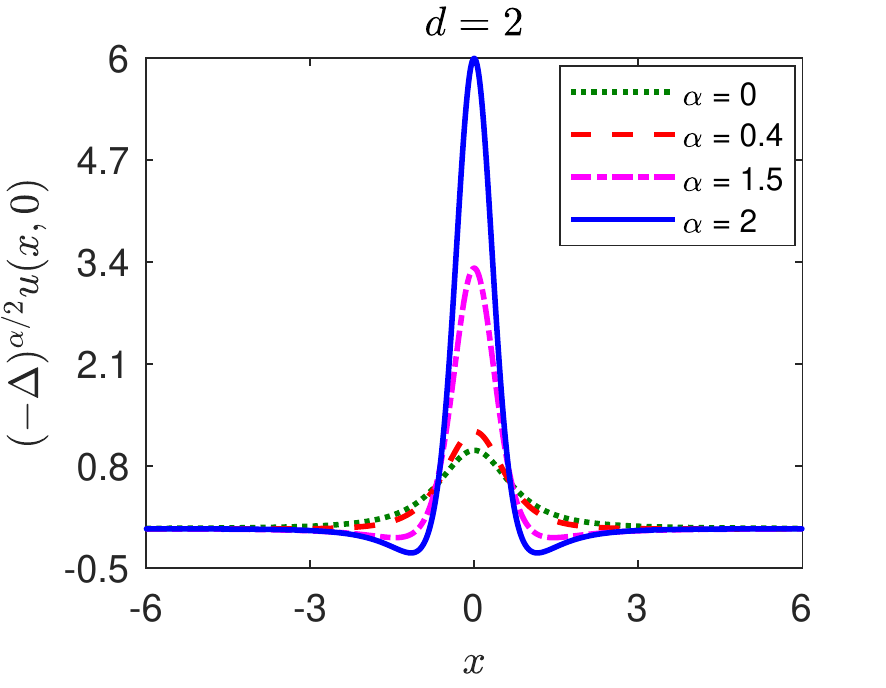} }
\caption{Illustration of $(-\Dt)^{\fl{\ap}{2}}u$ for the GIMQ  function $u(\bx) = \big(1+|\bx|^2\big)^{-(d+1)/2}$.}\label{Fig2-1}
\end{figure}

The discussion of GIMQ functions is recent, and thus their applications in solving PDEs remain very limited. 
It shows that for $\bt < -d/2$ the Fourier transform of GIMQ functions in \eqref{fun2-1-1}  can be given by the Mat\'ern function \cite{Wendland2005, Hamm2018}:
\bea\label{FGMQ}
\widehat{\varphi}(\xi) = \fl{2^{1+\bt} \veps^{2\bt}}{\Gamma(-\bt)}\big(\veps|\xi|\big)^{-(\bt+\fl{d}{2})}\, {K}_{-(\bt + \fl{d}{2})}\left(\fl{|\xi|}{\veps}\right), \qquad   \mbox{for} \ \ \xi \in {\mathbb R}^d\backslash\{{\bf 0}\},
\eea
where $K_v$ denotes the univariate modified Bessel function of the second kind, i.e., 
\beas
K_v(x) := \int_0^\infty e^{-x\, {\rm cosh(\tau)}} {\rm cosh}(v\tau) d\tau, \qquad  \mbox{for} \ \  v \in{\mathbb R}, \ \ x > 0.
\eeas
In this work, {\it we will consider  GIMQ functions \eqref{fun2-1-1} with a dimension-dependent power  $\bt = -(d+1)/2$}. 
In this case, the GIMQ function is also known as the Poisson kernel (up to a constant) \cite{Yoon2001}. 
Substituting  $\bt = -(d+1)/2$ into \eqref{FGMQ} reduces the Fourier transform as
\beas
\widehat{\varphi}(\xi) = \fl{2^{-d/2} \sqrt{\pi}}{\Gamma\big((d+1)/2\big)}e^{-|\xi|/\veps}, \qquad   \xi \in {\mathbb R}^d\backslash\{{\bf 0}\},
\eeas
by noticing  $K_{-\fl{1}{2}}(x) = K_\fl{1}{2}(x) = \sqrt{\pi/(2x)}\, e^{-x}$.
Next, we will introduce the Laplacian of GIMQ functions, which is an important building block of our method. 
\begin{lemma}[{\bf Laplacian of generalized inverse multiquadrics}] 
\label{lemma1}
Let the GIMQ function $u(\bx) = (1+|\bx|^2)^{-(d+1)/2}$, i.e. choosing $\bt = -(d+1)/2$  in (\ref{fun2-1-1}), for  dimension $d\ge 1$. 
Then, the function $(-\Dt)^{\fl{\ap}{2}}u$ can be analytically given by \cite{Dyda2012}: 
\bea\label{fun2-2-2}
(-\Dt)^{\fl{\ap}{2}} u(\bx) = \fl{2^{1-d}\sqrt{\pi} \Gamma\big(d+\ap\big)}{\Gamma(d/2)\Gamma((d+1)/2)} \,_2F_1\Big(\fl{d+\ap}{2}, \, \fl{d+1+\ap}{2};\, \fl{d}{2}; \, -|\bx|^2\Big), \qquad \mbox{for} \ \  \ap \ge 0, 
\eea
where $_2F_1$ denotes the Gauss hypergeometric function. 
\end{lemma}
Lemma \ref{lemma1} provides the analytical results of $(-\Dt)^{\fl{\ap}{2}} u$ for the GIMQ function $u(\bx) = (1+|\bx|^2)^{-(d+1)/2}$ at any point $\bx \in{\mathbb R}^d$. 
Moreover, this result holds for any power $\ap \ge 0$. 
In the special case of $\ap = 0$,  the right hand side of  \eqref{fun2-2-2} reduces to  the GIMQ function 
$u(\bx) = (1+|\bx|^2)^{-(d+1)/2}$. 
While $\ap = 2$,  it collapses to  
\beas
-\Dt u(\bx) = (d+1)\big(1+|\bx|^2\big)^{-\fl{d+5}{2}}\big[d - 3|\bx|^2\big], \qquad \mbox{for}\  \ \bx \in{\mathbb R}^d.
\eeas
Fig. \ref{Fig2-1} illustrates the results $(-\Dt)^\fl{\ap}{2}u$ for various $\ap$ and $d = 1, 2$. 
Note that when $d = 1$ the GIMQ function in Lemma \ref{lemma1} actually reduces to the inverse quadratic function. 
Fig. \ref{Fig2-1} shows that the larger the parameter  $\ap$, the faster the function $(-\Dt)^\fl{\ap}{2}u$ decays to zero as $|\bx| \to \infty$.  
For the same $\ap$,  the function $(-\Dt)^\fl{\ap}{2}u$ decays faster in higher dimension. 

The uniform representation of classical and fractional Laplacian of GIMQ functions in Lemma \ref{lemma1} 
 provides a foundation for developing compatible GIMQ-based schemes for these two operators. 
Moreover, the analytical formulation in (\ref{fun2-2-2}) allows us to avoid numerical approximation of the hypersingular integral for the fractional ($\ap < 2$) Laplacian. 
The Gaussian hypergeometric function $\,_2F_1$ can be  calculated with the well-established methods (see \cite{Pearson2017} and references therein). 
If dimension $d$ is odd (i.e., $d = 1$, or $3$), we can use the properties of the Gauss hypergeometric function $\,_2F_1$ and rewrite (\ref{fun2-2-2}) in terms of elementary functions. 
In one-dimensional ($d = 1$) case,  the result in \eqref{fun2-2-2} reduces to
\bea\label{1D1}
(-\Dt)^{\fl{\ap}{2}} \big(1+x^2\big)^{-1} = \Gamma(1+\ap)\cos\big[(1+\ap)\arctan|x|\big](1+x^2)^{-\fl{1+\ap}{2}}, \qquad \mbox{for}  \ \ x \in {\mathbb R}. 
\eea
Note that $\cos\big({\rm arctan}|x|\big) = \sqrt{1+x^2}$, and thus \eqref{1D1} collapses to function $u$ when $\ap = 0$.
While in three-dimensional ($d  = 3$) case, we have
\bea\label{3D3}
(-\Dt)^{\fl{\ap}{2}} \big(1+|\bx|^2\big)^{-2} = \left\{\begin{array}{ll}
\displaystyle\Gamma(3+\ap)/2, & \mbox{for} \ \ \bx = {\bf 0},\\ 
\displaystyle \fl{\Gamma(2+\ap)(1+|\bx|^2)^{-\fl{2+\ap}{2}}}{{2|\bx|}}{\sin\big((2+\ap)\arctan|\bx|\big)},\ \  & \mbox{otherwise},\\
\end{array}\right.
\eea
for $\bx \in {\mathbb R}^3$.  
Hence, when $d = 1$ or $3$,   one can use the results in (\ref{1D1}) and (\ref{3D3}) to replace  (\ref{fun2-2-2}) so as to avoid computing the Gauss hypergeometric function $\,_2F_1$.

Additionally, the Laplace operator satisfies the following properties \cite{Burkardt0020}: 
\bea
\label{prop1}
&&(-\Dt)^{\fl{\ap}{2}}\big[u(\bx - \bx_0)\big] = {\mathcal U}(\bx - \bx_0),  \quad\ \ \mbox{for} \ \ \bx_0\in{\mathbb R}^d, \qquad\qquad\qquad\\
\label{prop2}
&&(-\Dt)^{\fl{\ap}{2}}\big[u(\xi\bx)\big] = |\xi|^\ap{\mathcal U}(\xi\bx), \quad\ \ \mbox{for} \ \ \xi\in{\mathbb R}, 
\eea
for any $\ap \ge 0$,  where we denote function ${\mathcal U}(\bx) := (-\Dt)^{\fl{\ap}{2}}u(\bx)$. 
These two properties together with (\ref{fun2-2-2}) play an important role in the design of our meshless methods in Section \ref{section3}. 

\section{Meshless method with GIMQ}
\label{section3}
\setcounter{equation}{0}

In this section, we will present a new meshless method based on the GIMQ RBFs to solve the diffusion problem (\ref{diffusion})--(\ref{diffusion-BC}).  
Our method differs from other RBF-based methods in the following aspects. 
First, it provides an $\ap$-parametric scheme that can solve both classical  ($\ap = 2$) and fractional ($\ap <  2$) PDEs problems seamlessly. 
Second,  utilizing the properties of Laplace operators and GIMQ functions, our method avoids numerical evaluations of fractional derivatives, which are usually approximated by quadrature rules in many other methods \cite{Piret2013, Rosenfeld2019}. 
This not only reduces computational cost but  simplifies the implementations  especially in high dimensions. 
Last but not least, when $\ap < 2$ our method naturally integrates the extended boundary conditions into the scheme and avoids boundary truncation  as used in other method \cite{Rosenfeld2019}.

To introduce our method, let's first focus on the spatial discretization, i.e., approximating the Laplace operator $(-\Dt)^\fl{\ap}{2}$ for $\ap \in (0, 2]$.  
Assume that the solution of the diffusion problem (\ref{diffusion})--(\ref{diffusion-BC}) takes the ansatz: 
\bea\label{Sol1D}
u(\bx, t) \approx \widehat{u}(\bx, t) :=\sum_{i=1}^{\bar{N}} \lambda_i(t)\,\varphi^\veps(|\bx - \bx_i|), \qquad \mbox{for} \ \ \bx \in \bar{\Og}, \ \ t \ge 0, 
\eea
where $\bx_i$ are the RBF center points, and $\varphi^\veps(|\bx-\bx_i|)$ represents the GIMQ function centered at $\bx_i$. 
In our method, we choose the GIMQ basis function as
\beas
\varphi^\veps(|\bx - \bx_i|)= \big(1+\veps_i^2|\bx - \bx_i|^2\big)^{-\fl{d+1}{2}},
\eeas
where without loss of generality, we assume that the shape parameter $\veps_i$ associated with each center point $\bx_i$  is different. 
Note that the power $\bt = -(d+1)/2$, that is, the GIMQ basis functions change for different dimension $d \ge 1$. 
For all $\ap \in (0, 2]$, the RBF center points are assigned only on $\bar{\Og} = \Og \cup \p\Og$. 
More specifically, we assume that point $\bx_i \in \Og$ if index $1 \le i \le N$, while $\bx_i \in \p\Og$ if index $N+1 \le i \le \bar{N}$.

We start with approximating the fractional ($\ap < 2$) Laplacian subject to the Dirichlet boundary conditions in \eqref{diffusion-BC}. 
To this end, we consider the pointwise definition of the fractional Laplacian. 
Substituting the ansatz \eqref{Sol1D} into the fractional Laplacian (\ref{integralFL})  and taking the boundary conditions (\ref{diffusion-BC}) into account, we obtain
\bea\label{Dlaplace1}
&&(-\Dt)^{\fl{\ap}{2}}_hu(\bx, t) = C_{d, \ap} \bigg({\rm P.V.}\int_\Og\fl{\widehat{u}(\bx, t)- \widehat{u}(\by, t)}{|\bx - \by|^{d+\ap}} d\by + \int_{\Og^c}\fl{\widehat{u}(\bx, t) - g(\by, t)}{|\bx - \by|^{d+\ap}} d\by \bigg) \qquad\qquad\qquad \nn\\
&&\hspace{2.35cm}= (-\Dt)^\fl{\ap}{2} \widehat{u}(\bx, t) +  C_{d, \ap}  \int_{\Og^c}\fl{\widehat{u}(\by, t) - g(\by, t)}{|\bx - \by|^{d+\ap}} d\by, \qquad \mbox{for} \ \ \bx \in \Og. 
\eea
Note that the integral term over $\Og^c$ comes from the extended Dirichlet boundary conditions.  
Hence, it appears only in the fractional cases. 
On the other hand, we can directly apply the operator $-\Dt$ to  (\ref{Sol1D}) and obtain the approximation in the classical cases, i.e., $(-\Dt)_h u(\bx, t) = -\Dt \widehat{u}(\bx, t)$ for $\bx \in \Og$.  
Combining the classical and fractional cases,  we obtain a unified approximation: 
\bea\label{Eq-uniform}
(-\Dt)_h^{\fl{\ap}{2}}u(\bx, t) =  (-\Dt)^{\fl{\ap}{2}}\widehat{u}(\bx, t) +\zeta_\ap C_{d, \ap} \int_{\Og^c}\fl{\widehat{u}(\by, t) - g(\by, t)}{|\bx - \by|^{d+\ap}} d\by,\quad \ \mbox{for} \ \ \bx\in\Og, \ \ \ap \in (0, 2],
\eea
where $\zeta_\ap = 1 - \lfloor\ap/2\rfloor$ with $\lfloor \cdot \rfloor$ being the floor function. 
By Lemma \ref{lemma1} and properties \eqref{prop1}--\eqref{prop2}, it is easy to get
\bea\label{Eq-classical}
(-\Dt)^\fl{\ap}{2}\widehat{u}(\bx, t) = \underbrace{\fl{2^{1-d}\sqrt{\pi}\Gamma(d+\ap)}{\Gamma(d/2)\Gamma((d+1)/2)}}_{c_{d,\ap}}\sum_{i = 1}^{\bar{N}}\lambda_i(t)|\varepsilon_i|^\ap \,_2F_1\Big(\fl{d+\ap}{2}; \fl{d+1+\ap}{2}; \fl{d}{2};  -\varepsilon_i^\ap|\bx - \bx_i|^2\Big),
\eea
for $\bx \in \Og$ and $\ap \in (0, 2]$. 
Here, we denote the constant $c_{d, \ap}$ for notational simplicity, which is different from $C_{d, \ap}$ in the definition (\ref{integralFL}). 
It is evident that for $\ap \in (0, 2)$, the scheme \eqref{Eq-uniform}--\eqref{Eq-classical} exactly counts the extended Dirichlet boundary conditions via the integral term in (\ref{Eq-uniform}), and no extra unknowns are introduced on ${\mathbb R}^d\backslash\bar{\Og}$.

\begin{remark}
The scheme \eqref{Eq-uniform}--\eqref{Eq-classical}  can be straightforwardly applied to discretize the operator $(-\Dt)^{2m}$ with $m \in {\mathbb N}$. 
It can be also generalized to approximate the operator $(-\Dt)^{\fl{\ap}{2}}$ with $\ap > 2$ but $\ap \neq 2m$, but this generalization requires the point-wise definition of  $(-\Dt)^\fl{\ap}{2}$, such as the definition (\ref{integralFL}) for $\ap \in (0, 2)$, which is beyond the scope of this work. 
\end{remark}

Choose test points $\bx_k \in \bar{\Og}$ for $1 \le k \le \bar{N}$, that is, the total number of test points is the same as that of center points, but test points may not necessarily come from the same set of center points. 
More discussion on the choices of RBF center and test points can be found in \cite{Fornberg2002, Sarra2009, Fornberg2015}. 
Without loss of generality, we assume that test point $\bx_k \in \Og$ if index $1 \le k \le M$, while $\bx_k \in \p\Og$ if index $M+1 \le k \le \bar{N}$. 
For point $\bx_k \in \Og$,  the direct application of scheme  (\ref{Eq-uniform})--(\ref{Eq-classical}) to diffusion problem \eqref{diffusion} yields the semi-discretization form as: 
\bea\label{scheme1}
&&\sum_{i = 1}^{\bar{N}}\fl{d\lambda_i(t)}{dt} \varphi^{\veps}(|\bx_k - \bx_i|^2) = -\kappa\sum_{i = 1}^{\bar{N}} \lambda_i(t)\bigg[c_{d, \ap}|\varepsilon_i|^\ap\,_2F_1\Big(\fl{d+\ap}{2}, \, \fl{d+\ap+1}{2};\, \fl{d}{2}; \, -\veps_i^2|\bx_k - \bx_i|^2\Big)\quad   \nn\\
&&\hspace{1.5cm} +\zeta_\ap C_{d, \ap}\int_{\Og^c}\fl{\varphi^\veps(|\by-\bx_i|)}{|\bx_k - \by|^{d+\ap}} d\by\bigg] + \zeta_\ap C_{d,\ap}\int_{\Og^c}\fl{g(\by, t)}{|\bx_k-\by|^{d+\ap}}d\by+f\big(\bx_k, t \big), \quad  t > 0, \quad
\eea
for $k = 1, 2, \ldots, M$. 
While for test point $\bx_k \in \p\Og$, we can  combine \eqref{diffusion-BC} and (\ref{Sol1D}) to obtain the discretization of boundary conditions as: 
\bea\label{scheme2}
\sum_{i = 1}^{\bar{N}}\lambda_i(t) \varphi^{\veps}(|\bx_k - \bx_i|^2) = g(\bx_k, t), \quad  t \ge 0.
\eea
for $k = M+1, M+2, \ldots, \bar{N}$. 
It shows that boundary discretization (\ref{scheme2}) is only carried out for points on $\p\Og$ (instead of on $\Og^c$),  same for both  classical and fractional cases. 
If $\ap < 2$, the boundary conditions in (\ref{diffusion-BC}) for $\bx \in {\mathbb R}^d\backslash\bar{\Og}$ have been already applied to  scheme (\ref{scheme1}) via its last integral.  
The initial condition at time $t = 0$ can be discretized as:
\bea\label{scheme3}
\sum_{i = 1}^{\bar{N}}\lambda_i(0)\,\varphi^{\veps}(|\bx_k - \bx_i|^2) = u(\bx_k, 0),\qquad\mbox{for $k = 1,\, 2,\, \ldots,\, \bar{N}$.}
\eea
The scheme in \eqref{scheme1}--\eqref{scheme3} provides an ODE system for unknowns $\lambda_i(t)$ for $1 \le i \le \bar{N}$, which can be solved by any time stepping method. 

In the following, we will discretize the ODE system \eqref{scheme1}--\eqref{scheme3} by the Crank--Nicolson method. 
We first denote vector ${\boldsymbol \lambda}(t) = \big(\lambda_1(t), \, \lambda_2(t), \, \cdots, \,  \lambda_{\bar{N}}(t)\big)^T$ and rewrite system \eqref{scheme1}--\eqref{scheme2}  into matrix-vector form:
\bea\label{scheme4}
\Phi_{(1:M,\,1: \bar{N})}\fl{d{\boldsymbol \lambda}(t)}{dt} = -\kappa A_{M \times \bar{N}}\,{\boldsymbol \lambda}(t) + {\bf b}(t); \qquad \Phi_{(M+1: \bar{N},\,1: \bar{N})}\,{\boldsymbol \lambda}(t) = {\bf g}(t). 
\eea
Here,  $\Phi = \big\{\Phi_{k, i}\big\}_{\bar{N}\times \bar{N}}$ is a matrix of GIMQ basis functions with entry $\Phi_{k, i} = \varphi^\veps\big(|\bx_k - \bx_i|\big)$, and $\Phi_{(r_1 : r_2, \, c_1 : c_2)}$ denotes its submatrix including row $r_1$ to $r_2$ and column $c_1$  to $c_2$  from matrix $\Phi$. 
While $A$ represents the differentiation matrix of $(-\Dt)^{\fl{\ap}{2}}$ with its entry
\beas
A_{k, i} = c_{d, \ap}|\varepsilon_i|^\ap\,_2F_1\Big(\fl{d+\ap}{2}, \, \fl{d+\ap+1}{2};\, \fl{d}{2}; \, -\veps_i^2|\bx_k - \bx_i|^2\Big) + \zeta_\ap C_{d, \ap}\int_{\Og^c}\fl{\varphi^\veps(|\by-\bx_i|)}{|\bx_k - \by|^{d+\ap}} d\by,
\eeas
for $1 \le k \le M$ and $1 \le i \le \bar{N}$. 
The vector ${\bf g}_{(\bar{N} - M)\times 1} = \big(g(\bx_{M+1},  t),\,g(\bx_{M+2},  t), \,\ldots, \, g(\bx_{\bar{N}},  t)\big)^T$, and vector ${\bf b}_{M \times 1}$ has entry 
\beas
b_k(t) = f(\bx_k, t) +  \zeta_\ap C_{d,\ap}\int_{\Og^c}\fl{g(\by, t)}{|\bx_k-\by|^{d+\ap}}d\by, \quad \ \  \mbox{for} \ \ k = 1, 2, \ldots, M.
\eeas
Denote time sequence $t_n = n\tau$ (for $n = 0, 1, \ldots$) with time step $\tau > 0$. 
Then the Crank--Nicolson discretization of \eqref{scheme4} yield the fully discrete scheme as 
\bea\label{scheme5}
\left(\begin{array}{c}\displaystyle \Phi_{(1:M,\,1:\bar{N})} + \fl{\kappa\tau}{2} A\\
\displaystyle \Phi_{(M+1:\bar{N},\,1:\bar{N})}\end{array}\right)
{\boldsymbol \lambda}^{n+1} = 
\left(\begin{array}{c}\displaystyle \Big(\Phi_{(1:M,\,1:\bar{N})} - \fl{\kappa\tau}{2} A\Big){\boldsymbol\lambda}^n+ \fl{\tau}{2}\big({\bf b}(t_n) + {\bf b}(t_{n+1})\big) \\
\displaystyle{\bf g}(t_{n+1})\end{array}\right), 
\eea
for $n = 0, 1, \ldots$, where ${\boldsymbol \lambda}^n$ represents the numerical approximation of ${\boldsymbol \lambda}(t_n)$. 
Initially, we can obtain ${\boldsymbol \lambda}^0$ by solving $\Phi_{\bar{N}\times\bar{N}} {\boldsymbol \lambda}^0 = {\bf U}_0$ with ${\bf U}_0 = \big(u(\bx_1, 0), \, u(\bx_2, 0), \, \ldots, \, u(\bx_{\bar{N}}, 0)\big)^T$. 
Substituting ${\boldsymbol \lambda}^n$ into \eqref{Sol1D} leads to the numerical solution of \eqref{diffusion}--\eqref{diffusion-BC} at time $t_n$ and for any point $\bx \in {\Og}$.

Since GIMQ functions are globally supported, the linear system (\ref{scheme5}) has a full  stiffness matrix  for both classical ($\ap = 2$) and fractional ($\ap < 2$) problems. 
On the other hand, the nonlocality of the fractional Laplacian always leads to a linear system with full matrix even if it is discretized by local (e.g., finite difference/element) methods \cite{Duo2018, Duo-FDM2019, Acosta2017}. 
Hence,  our method does not introduce extra computations in the fractional cases. 
Instead, it could save computational cost by using fewer number of points  to achieve the desired accuracy. 
This suggests that global methods might be more beneficial for solving nonlocal or fractional problems.

\section{Accuracy, comparison, and shape parameters}
\label{section4}
\setcounter{equation}{0}

In this section, we will test the performance of our  method  and compare it with the Gaussian RBF-based method recently proposed in the literature \cite{Burkardt0020}. 
Furthermore, we will study two approaches of selecting the shape parameter $\veps$ for better performance of our method. 
Here, we will focus on the one-dimensional Poisson problem: 
\bea
\label{Poisson}
(-\Dt)^{\fl{\ap}{2}} u = f(x), \quad \mbox{for} \ \,  x \in \Og; \qquad\ u(x) = g(x), &&  \mbox{for} \ \, x \in \Upsilon.
\eea
Choose the domain $\Og = (-1, 1)$. 
The definition of $\Upsilon$ implies that if $\ap = 2$ a two-point Dirichlet boundary condition is imposed at $x = \pm 1$, while an extended Dirichlet boundary condition is applied on $\Og^c = (-\infty, -1] \cup [1, \infty)$ if $\ap < 2$. 
For the purpose of easy comparison,  we will choose the RBF center points uniformly distributed on $[-1, 1]$ and  test points from the same set of center points (consequently, $M = N$). 
Note that even though the extended Dirichlet boundary conditions are imposed on $\Og^c$ in the fractional cases, the RBF points are only considered on $\bar{\Og}$. 
We remark that our method is flexible in choosing the center and test points, and the above choice is only an example that we use here. 
Unless otherwise stated, we will use a constant shape parameter, i.e., $\veps_i \equiv \veps$ for $1 \le i \le \bar{N}$,  in the following simulations.

Let $u_j$ and $u_j^h$ represent the exact solution and numerical approximation of $u(\bx)$ at point $\bx = \bx_j$, respectively.
In the following, we compute numerical errors  as the root mean square (RMS)  errors:
\beas\label{rms}
\|e\|_{\rm rms} = \bigg(\fl{1}{K}\sum_{j = 1}^K |u_j - u_j^h|^2\bigg)^{1/2},
\eeas
where  $K \gg \bar{N}$ denotes the total number of interpolation points on domain $\bar{\Og}$. 
In practice, a sufficiently large number $K$ is chosen such that the numerical error $\|e\|_{\rm rms}$ is insensitive to it. 
As we will see later, our method not only unifies the numerical discretization of the fractional  ($0<\alpha<2$) and classical ($\alpha=2$)  Laplacian in a single scheme, but also allows simple computer implementation  for any dimension $d \ge 1$. 

\subsection{Numerical accuracy}
\label{section4-1}

We will test the performance of our method in solving the Poisson problem \eqref{Poisson} with different boundary conditions. 
Most of the existing studies on fractional PDE problems focus on homogeneous boundary conditions (e.g. $g(x) \equiv 0$ in \eqref{Poisson}), while the results on nonhomogeneous boundary conditions remain very limited. 
One of the main challenges is to accurately count the nonzero boundary conditions into the scheme. 

\subsubsection{Homogeneous boundary conditions}
\label{section4-1-1}

Here, we consider a benchmark (fractional) Poisson problem with homogeneous boundary conditions, i.e., $g(x) \equiv 0$, and choose function $f$ in (\ref{Poisson}) as
\beas
f(x) = \frac{2^\alpha \Gamma\Big(\frac{\alpha+1}{2}\Big) \Gamma(s+1+\frac{\alpha}{2})} {\sqrt{\pi}\Gamma(s+1)}\,_2F_1\Big(\frac{\alpha+1}{2}, -s, \frac{1}{2}, x^2\Big), && \mbox{for} \ \,  s\in{\mathbb N}^0.
\eeas
For any $\ap \in (0, 2]$, the exact solution of this Poisson problem is given by $u(x) = (1-x^2)^{s+\fl{\alpha}{2}}_+$, a compact support function on $[-1, 1]$.
In this example, we will take $s = 3$.

Table \ref{Tab4-1-1} demonstrates the RMS errors and condition numbers $\mathcal{K}$ of the resultant linear system, where a constant shape parameter $\veps$ is used for all RBF center points. 
It shows that the numerical errors decrease with a spectral rate as  the number of points $\bar{N}$ increases. 
Our method can achieve a good accuracy even with fewer points.  
For instance, it has an error of ${\mathcal O}(10^{-7}) \sim {\mathcal O}(10^{-8})$ for $\bar{N} = 65$, which is  much smaller than the errors from finite difference/element methods  in the literature \cite{Duo2018,  Acosta2017}. 
\begin{table}[htb!]
\begin{center} 
\begin{tabular}{|c|c|c|c|c|c|c|c|c|}
\hline
& \multicolumn{2}{|c|}{$\ap = 0.6$, $\varepsilon = 3$} & \multicolumn{2}{|c|}{$\ap = 1$, $\varepsilon = 3.5$} & \multicolumn{2}{|c|}{$\ap = 1.5$, $\varepsilon = 3.5$} & \multicolumn{2}{|c|}{$\ap = 2$, $\varepsilon = 3.5$}   \\
\cline{2-9}
$\bar{N}$     & $\|e\|_{\rm rms}$ & $\mathcal{K}$ & $\|e\|_{\rm rms}$ & $\mathcal{K}$ & $\|e\|_{\rm rms}$ & $\mathcal{K}$ & $\|e\|_{\rm rms}$ & $\mathcal{K}$ \\ 
\hline
5& 1.233e-1 & 2.099 & 2.650e-1 & 4.086 & 3.838e-1 & 1.05e01 & 4.626e-1 & 2.87e01\\
\hline
9& 3.608e-3 & 6.542 & 2.616e-2 & 5.808 & 8.189e-2 & 1.07e01 & 1.980e-1 & 2.84e01\\
\hline
17 &2.468e-4 & 2.24e02 & 4.125e-4 & 4.92e01 & 8.420e-4 & 5.34e01 & 2.180e-3 & 5.24e01 \\
\hline
33& 1.866e-5 & 6.24e05 & 3.375e-5 & 2.25e04 & 5.933e-5 & 7.37e03 & 1.307e-4 & 2.03e03\\
\hline
65 & 8.655e-8 & 7.66e12 & 5.891e-8 & 2.16e10 & 3.355e-7 & 5.98e09 & 3.513e-7 & 8.47e08\\
\hline
\end{tabular}
\caption{Numerical errors  $\|e\|_{\rm rms}$ and condition numbers $\mathcal{K}$ in solving the 1D Poisson problem  (\ref{Poisson}), where exact solution is $u(x) = (1-x^2)^{3+\fl{\ap}{2}}_+$.}\label{Tab4-1-1}
\end{center}
\end{table}
Moreover,  our method has a distinct merit -- unifying the approximation scheme for the fractional and classical Laplacians in a single scheme. 
Additionally, we find that as $\bar{N}$ increases, the condition number of the differentiation matrix increases quickly, which indicates that the system could become ill-conditioned when a large number of points is used. 
Different strategies can be found in the literature to resolve the ill-conditioning issues of RBF-related methods; see \cite{Sarra2011, Sarra2009, Fornberg2015} and references therein. 
In Section \ref{section4-3}, we will study two approaches to suppress the condition number via shape parameter $\veps$. 

\begin{figure}[htb!]
\centerline{\includegraphics[height = 5.76cm, width = 7.86cm]{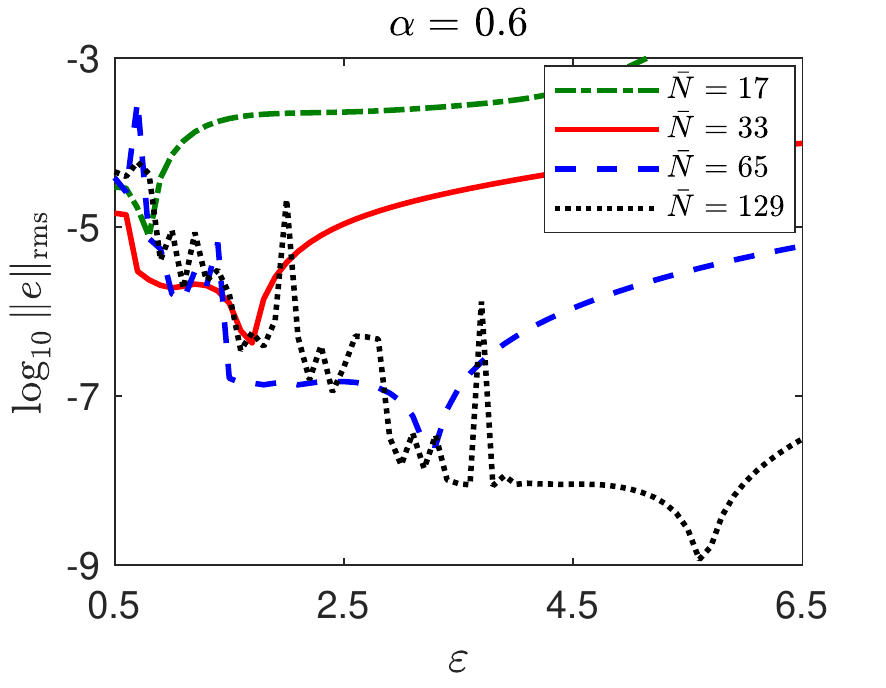}\hspace{-5mm}
\includegraphics[height = 5.76cm, width = 7.36cm]{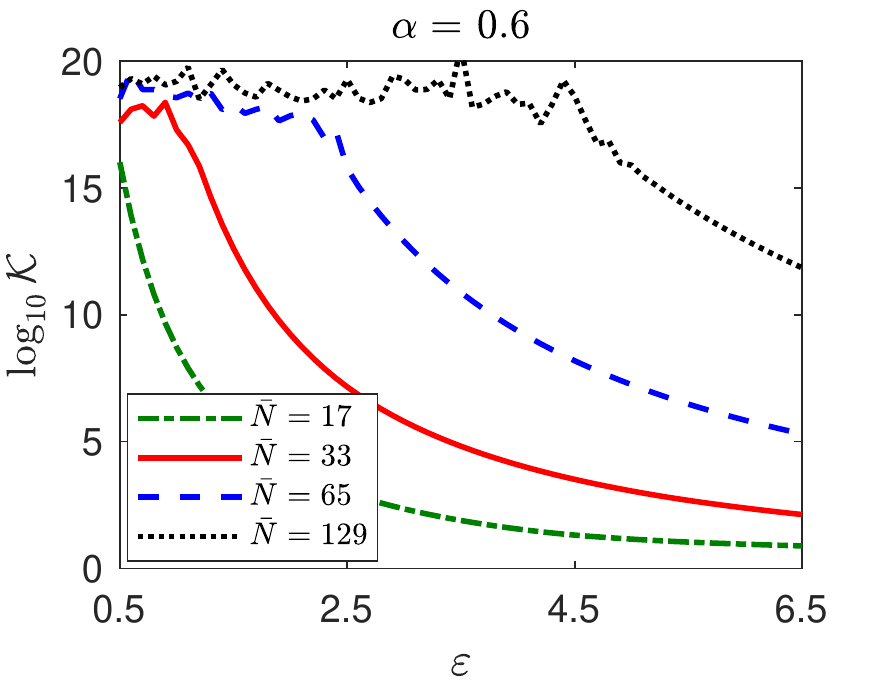}}
\centerline{\includegraphics[height = 5.76cm, width = 7.86cm]{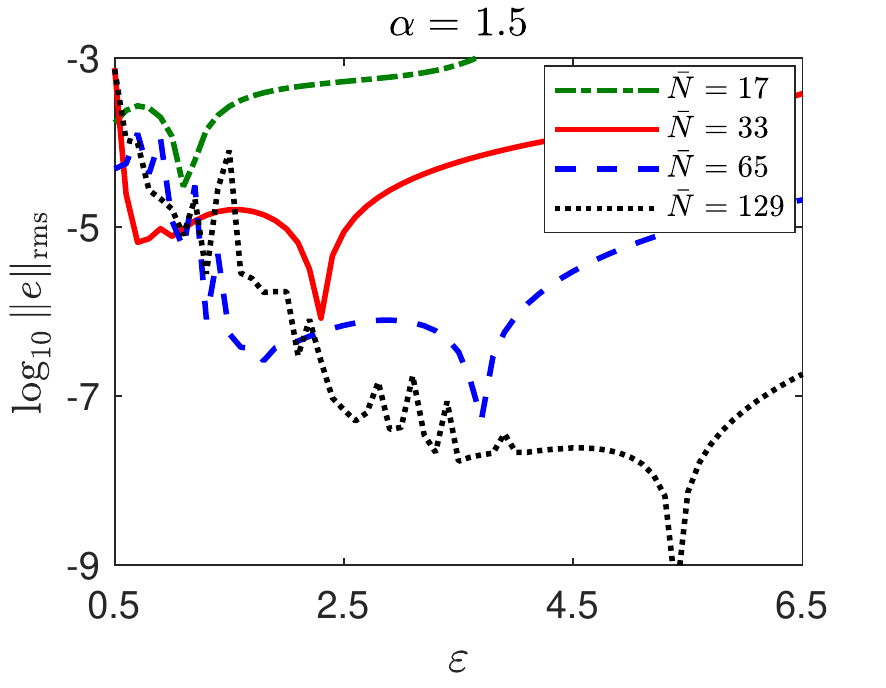}\hspace{-5mm}
\includegraphics[height = 5.76cm, width = 7.36cm]{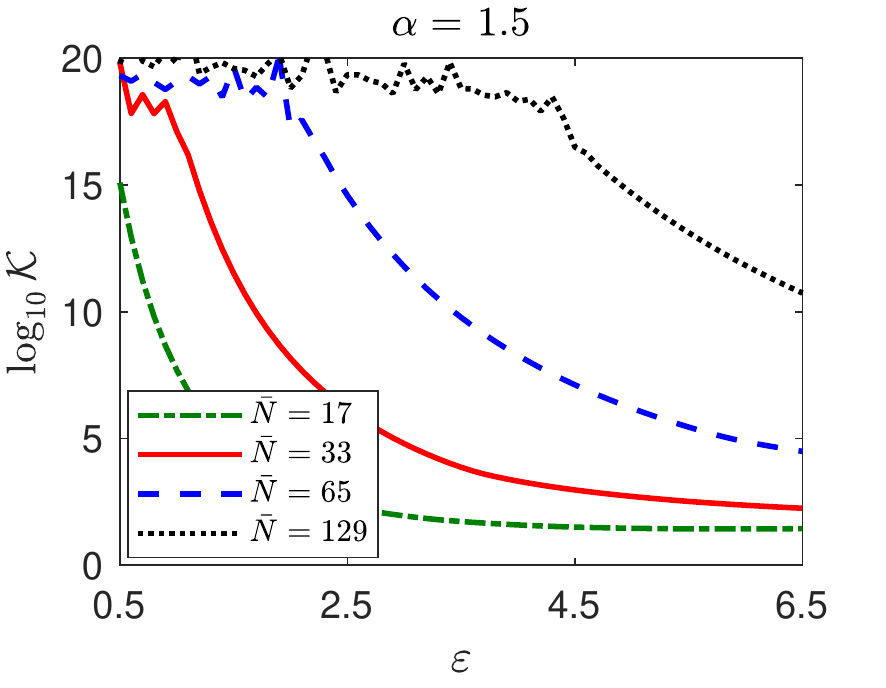}}
\caption{Numerical errors $\|e\|_{\rm rms}$ and condition numbers $\mathcal{K}$ versus the constant shape parameter $\veps$ in solving the 1D Poisson problem (\ref{Poisson}) with  exact solution $u(x) = (1-x^2)^{3+\fl{\ap}{2}}_+$.} \label{Fig3-1-1}
\end{figure}

In Fig. \ref{Fig3-1-1}, we further study the relation between numerical errors $\|e\|_{\rm rms}$, condition number ${\mathcal K}$, total number of points $\bar{N}$, and the constant shape parameter $\veps$.  
We find that for a given $\bar{N}$,   numerical errors generally decrease first and then increase with the shape parameter $\veps$. 
There exists an optimal shape parameter $\veps^*({\bar{N}})$,  depending on the total number of points and solution property, at which the numerical error is minimized. 
The condition number could be very large if the number of points $\bar{N}$ is large or the shape parameter $\veps$ is small. 
If  $\veps$ is large enough,  the condition number would monotonically decrease as $\veps$ increases. 
This implies that the ill-conditioning issue caused by a large number of points could be relaxed by increasing the shape parameter, but a large shape parameter might reduce numerical accuracy (see Fig. \ref{Fig3-1-1} left panel). 

\subsubsection{Nonhomogeneous boundary conditions}
\label{section4-1-2}

In the literature, the results on fractional PDEs with nonhomogeneous Dirichlet boundary conditions still remain rare.  
In this example, we consider the Poisson problem (\ref{Poisson}) with 
\beas
f(x) = \frac{\sqrt{2}}{(1+\ap)\sqrt{\pi}}\,_2F_1\Big(\frac{\alpha+1}{2}; \frac{3+\ap}{2}, \frac{1}{2}; -\frac{1}{4}x^2\Big) , \qquad g(x) = \sqrt{\fl{2}{\pi}}{\rm sinc}\Big(\fl{x}{\pi}\Big).
\eeas
Its exact solution is given by $u(x) = \sqrt{2/\pi}\,{\rm sinc}(x/\pi)$ for any $x \in {\mathbb R}$ and $\ap \in (0, 2]$. 
Table \ref{Tab4-1-2} shows the numerical errors and condition numbers for various $\ap$.

\begin{table}[htb!]
\begin{center} 
\begin{tabular}{|c|c|c|c|c|c|c|c|c|}
\hline
& \multicolumn{2}{|c|}{$\ap = 0.6$, $\varepsilon = 1$} & \multicolumn{2}{|c|}{$\ap = 1$, $\varepsilon = 1$} & \multicolumn{2}{|c|}{$\ap = 1.5$, $\varepsilon = 1.5$} & \multicolumn{2}{|c|}{$\ap = 2$, $\varepsilon = 1.5$}   \\
\cline{2-9}
$\bar{N}$     & $\|e\|_{\rm rms}$ & $\mathcal{K}$ & $\|e\|_{\rm rms}$ & $\mathcal{K}$ & $\|e\|_{\rm rms}$ & $\mathcal{K}$ & $\|e\|_{\rm rms}$ & $\mathcal{K}$ \\ 
\hline
5   &2.794e-3 &2.89e01 &3.506e-3 &1.64e01 &8.423e-3 &3.933  &9.620e-3  &5.678 \\
\hline
9   &2.355e-4 &5.50e03 &3.248e-4 &2.31e03 &3.387e-3  &3.82e01 &6.131e-3  &1.83e01\\
\hline
17 &2.855e-6 &5.13e08 &4.709e-6 &1.77e08 &2.366e-4  &4.20e04 &5.757e-4  &9.74e03\\
\hline
33 &1.360e-9 &1.89e17 &2.266e-9 &3.75e17 &1.907e-6  &3.35e11 &6.336e-6  &3.36e10\\
\hline
65 &6.406e-10 &3.51e18 &7.852e-10 &2.28e18  &9.585e-9 &4.65e19 &6.858e-9 &1.01e18\\
\hline
\end{tabular}
\caption{Numerical errors  $\|e\|_{\rm rms}$ and condition numbers $\mathcal{K}$ in solving the 1D Poisson problem  (\ref{Poisson}), where exact solution is $u(x) = \sqrt{2/\pi}\,{\rm sinc}(x)$ for $x \in {\mathbb R}$.}\label{Tab4-1-2}
\end{center}
\end{table}
We find similar observations as in Table \ref{Tab4-1-1} -- numerical errors decrease with a spectral rate as $\bar{N}$ increases. 
But,  the condition numbers in Table \ref{Tab4-1-2} are generally larger because of smaller shape parameters $\veps$ are used here. 
The computational cost in this example is higher due to extra efforts to evaluate the integrals of nonzero boundary conditions.  
\begin{figure}[htb!]
\centerline{
\includegraphics[height = 5.76cm, width = 7.86cm]{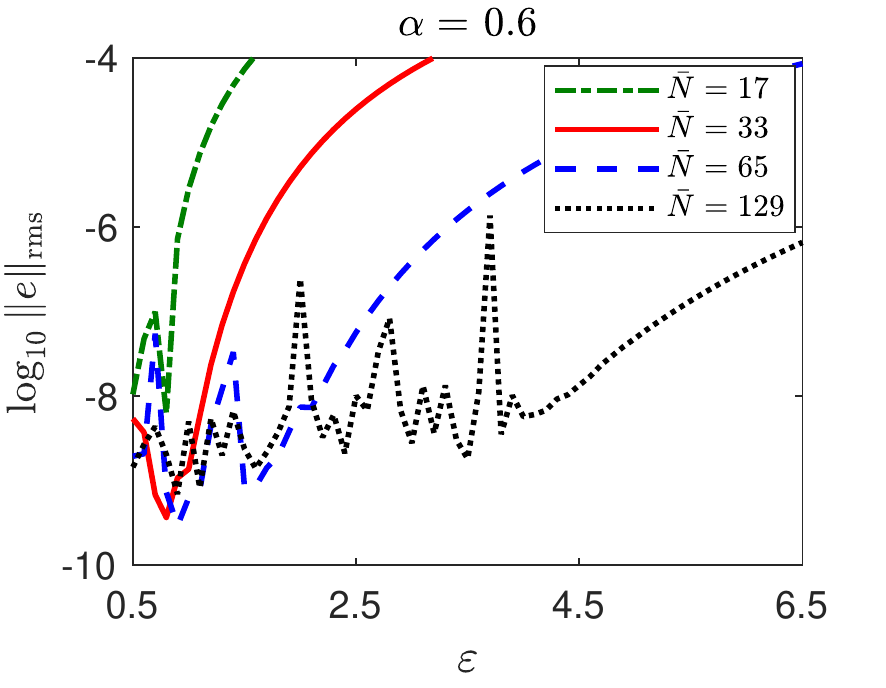}\hspace{-5mm}
\includegraphics[height = 5.76cm, width = 7.86cm]{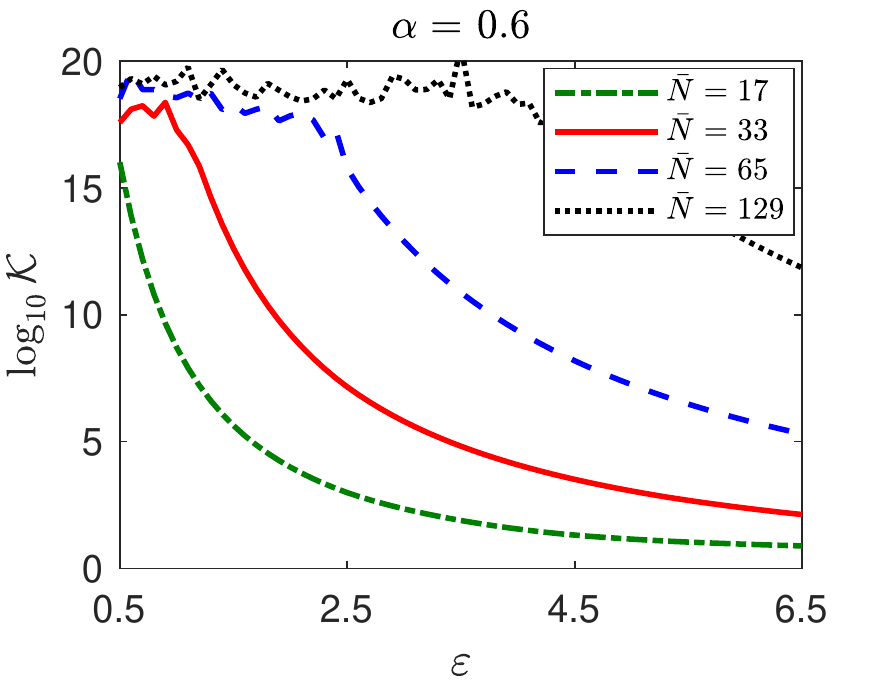}}
\centerline{
\includegraphics[height = 5.76cm, width = 7.86cm]{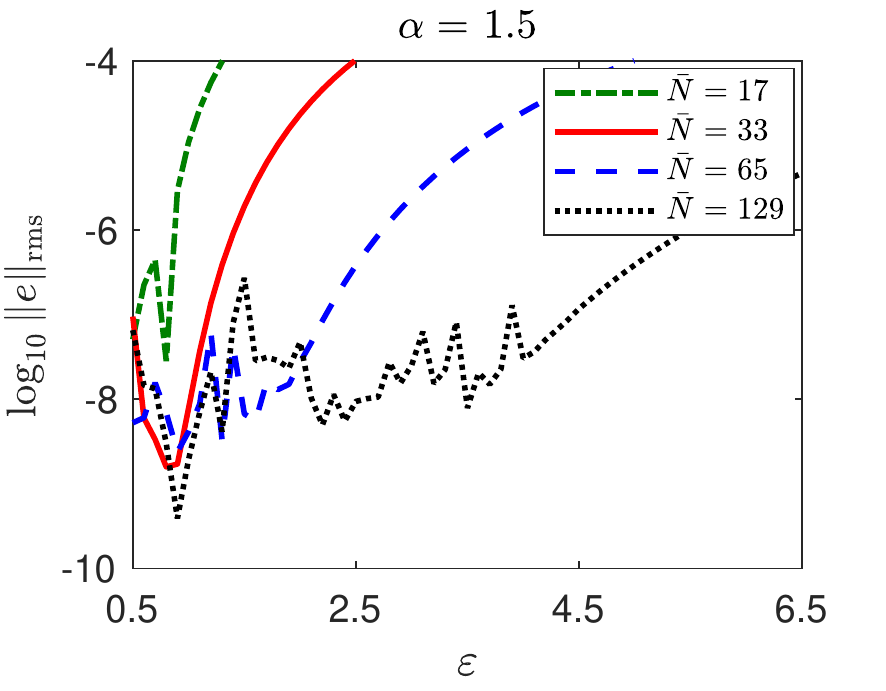}\hspace{-5mm}
\includegraphics[height = 5.76cm, width = 7.86cm]{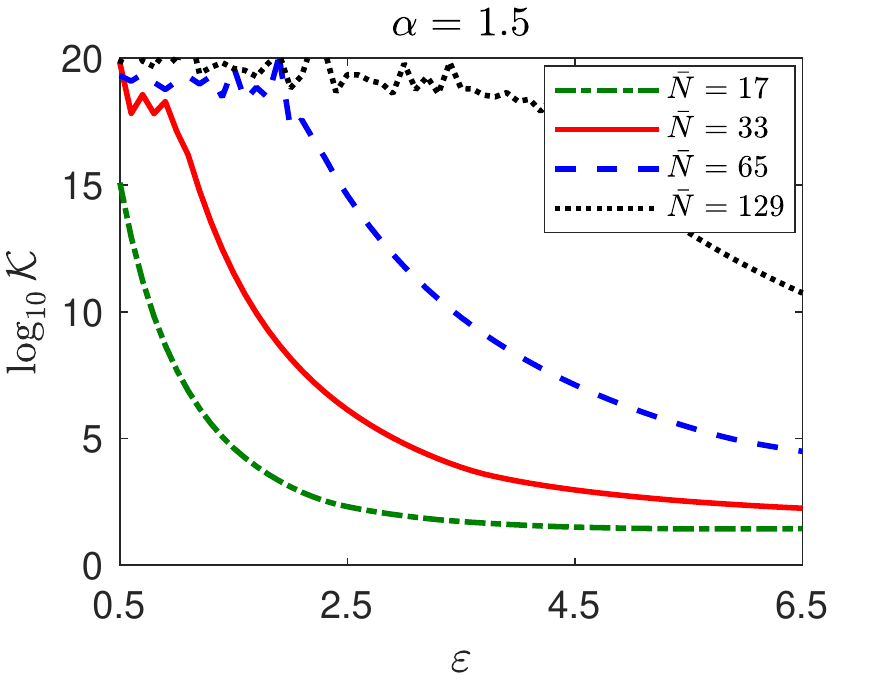}}
\caption{Numerical errors $\|e\|_{\rm rms}$ and condition numbers $\mathcal{K}$ versus the shape parameter $\veps$ in solving the 1D Poisson problem (\ref{Poisson}) with  exact solution $u(x) = \sqrt{2/\pi}\,{\rm sinc}(x)$ for $x \in {\mathbb R}$.} \label{Fig3-1-2}
\end{figure}
Fig. \ref{Fig3-1-2} demonstrates the relation of numerical errors $\|e\|_{\rm rms}$, conditional number ${\mathcal K}$, total number of points $\bar{N}$, and the shape parameter $\veps$.  
Similar to observations in Fig. \ref{Fig3-1-1}, numerical errors first decrease and then increase with the shape parameter $\veps$, but the optimal shape parameter in this example is less sensitive to the number of points $\bar{N}$. 
The minimum errors are achieved at a shape parameter $\veps^*(\bar{N})$ around $0.5 \sim 1.5$ for all $\bar{N}$ and $\ap$  (see Fig. \ref{Fig3-1-2} left panel). 
Notice that the exact solution $u \in C^{3, \fl{\ap}{2}}({\bar{\Og}}$) in Fig. \ref{Fig3-1-1}, while $u \in C^{\infty}(\bar{\Og})$ in Fig. \ref{Fig3-1-2}, suggesting that solution regularity might play an important role in determining the optimal shape parameters.  
The above studies also suggest that it is challenging to find a uniform optimal shape parameter $\veps$ for different numbers of points $\bar{N}$. 
In practice, it is important to select a shape parameter that can yield good numerical accuracy and avoid ill-conditioning issues in the computation; see more discussion in Section \ref{section4-3}. 

\subsection{Comparison to the Gaussian-based method}
\label{section4-2} 

Recently, a Gaussian-based method was proposed in \cite{Burkardt0020} to discretize the classical and fractional Laplacians. 
As noted previously, the Gaussian and GIMQ functions share many similarities, i.e.,  globally supported, infinitely differentiable, and strictly positive definite. 
Both of them have been well studied in interpolation problems \cite{Franke1979, Fornberg2008}. 
However, very few comparison of these two functions can be found in the literature, especially for solving PDEs. 
In the following, we will compare our GIMQ-based method with the Gaussian-based method in \cite{Burkardt0020} by studying their numerical errors and condition numbers with respect to constant shape parameters.

\begin{figure}[htb!]
\centerline{\includegraphics[height = 5.76cm, width = 7.86cm]{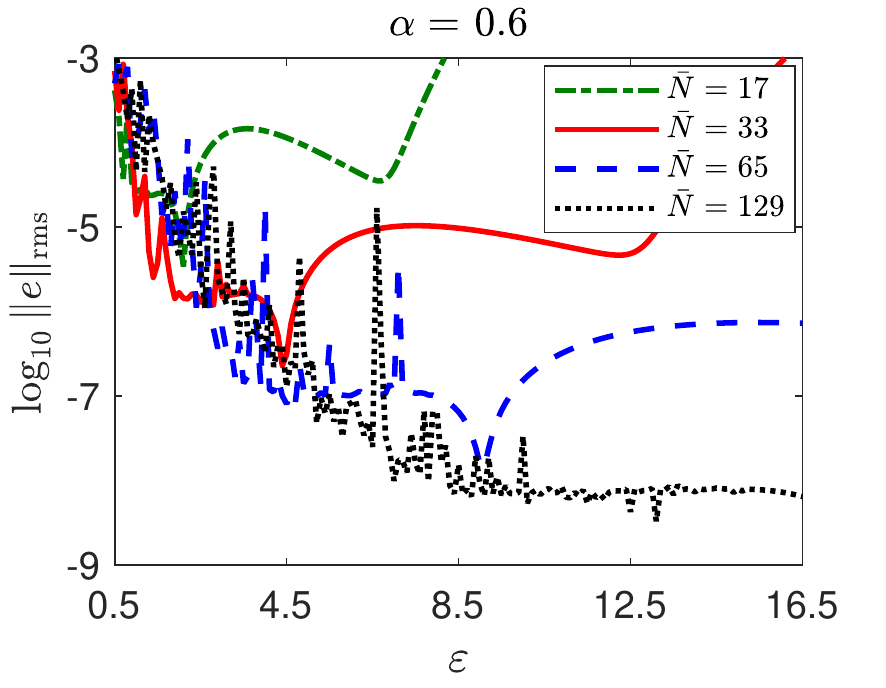}\hspace{-5mm}
\includegraphics[height = 5.76cm, width = 7.36cm]{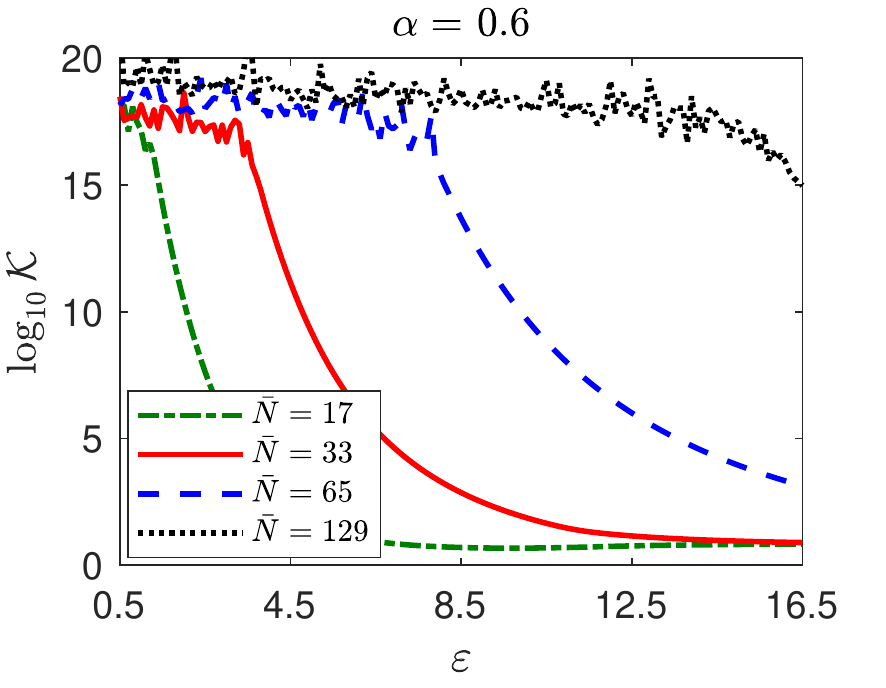}}
\centerline{\includegraphics[height = 5.76cm, width = 7.86cm]{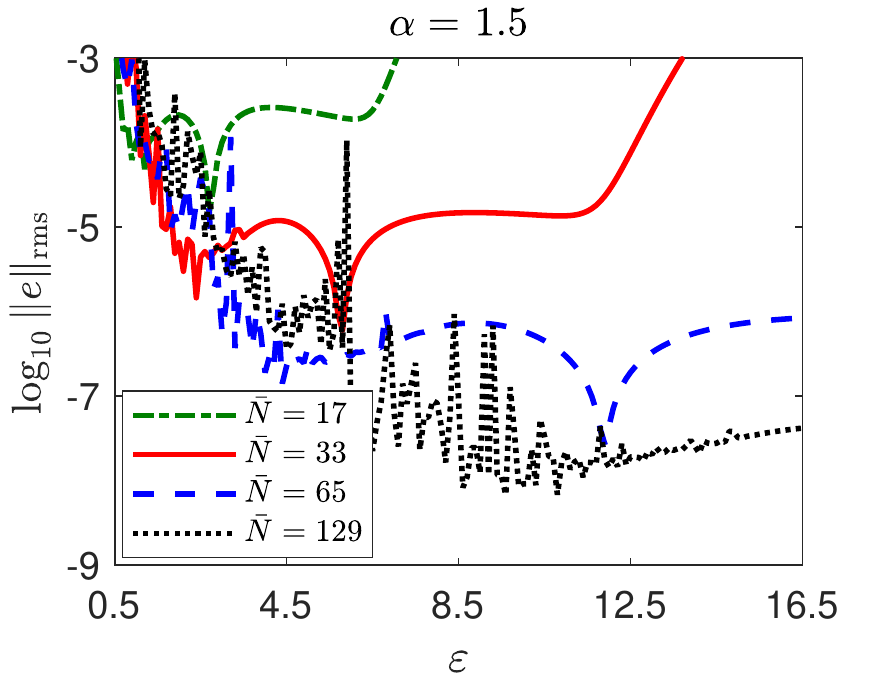}\hspace{-5mm}
\includegraphics[height = 5.76cm, width = 7.36cm]{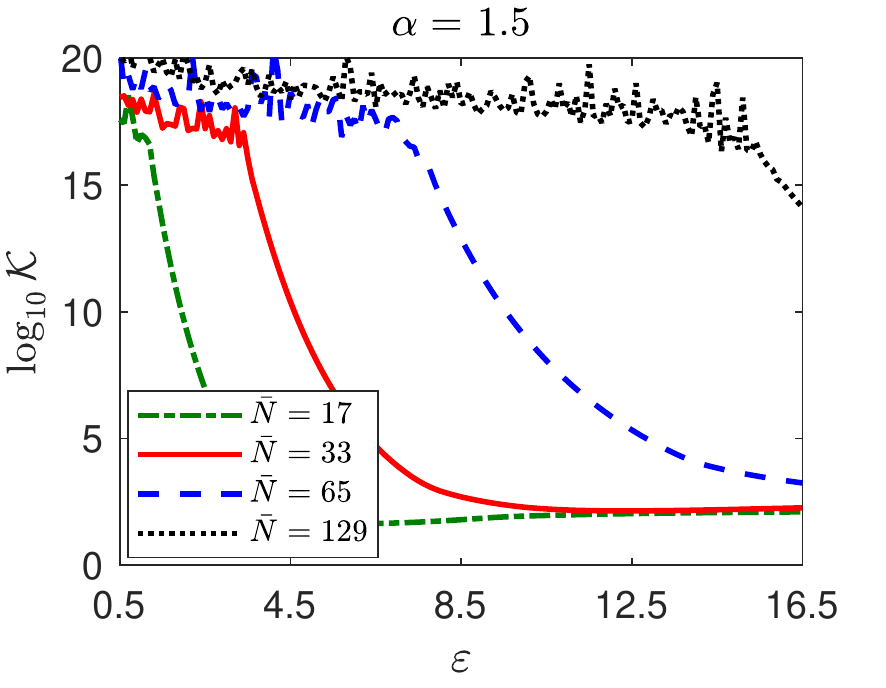}}
\caption{Gaussian RBF-based method:  Numerical errors $\|e\|_{\rm rms}$ and condition numbers $\mathcal{K}$ versus the shape parameter $\veps$ in solving the 1D Poisson problem with  exact solution $u(x) = (1-x^2)^{3+\fl{\ap}{2}}_+$. These results are in comparison with those in Fig. \ref{Fig3-1-1} for the GIMQ-based method.  }\label{Fig4-2-1}
\end{figure}

Here, we will mainly focus on the Poisson problems as studied in Section \ref{section4-1-1}, and the constant shape parameters are used for all RBF center points. 
Our studies show that the Gaussian RBF-based method has a spectral accuracy. 
To compare with GIMQ results in Fig. \ref{Fig3-1-1}, we present in Fig. \ref{Fig4-2-1} the relation between numerical errors, condition number, number of points, and shape parameters of the Gaussian RBF-based method. 
It shows that  numerical errors of both methods first decrease and then increase with the constant shape parameter $\veps$, but the errors of GIMQ method decrease with less oscillation. 
There also exists an optimal constant shape parameter $\veps^* = \veps^*(\bar{N})$ for the Gaussian-based method. 
These two methods are different mainly in two aspects: 
(i) under the same numerical setting,  the optimal shape parameter $\veps^*$ of the Gaussian method is much larger.
Also, it is more sensitive to the number of points $\bar{N}$. 
(ii) The condition numbers of both methods monotonically decrease after a large shape parameter (see the  right panels of Figs. 3 and 5).  
But, the GIMQ-based method starts this monotonic decay at a smaller shape parameter, but the condition number of the Gaussian-based method remain large for a wide range of $\veps$. 
For instance, when $\ap = 0.6$ and $\bar{N} = 129$, the condition number reduces to ${\mathcal K} \sim 10^{15}$ at $\veps = 5$ in GIMQ-based method (see Fig. \ref{Fig3-1-1}), while $\veps = 16.5$ in Gaussian-based method (see Fig.  \ref{Fig4-2-1}).  
This suggests that it might be easier to suppress ill-conditioning issues in the GIMQ-based method.

\begin{figure}[htb!]
\centerline{
\includegraphics[height = 5.76cm, width = 7.86cm]{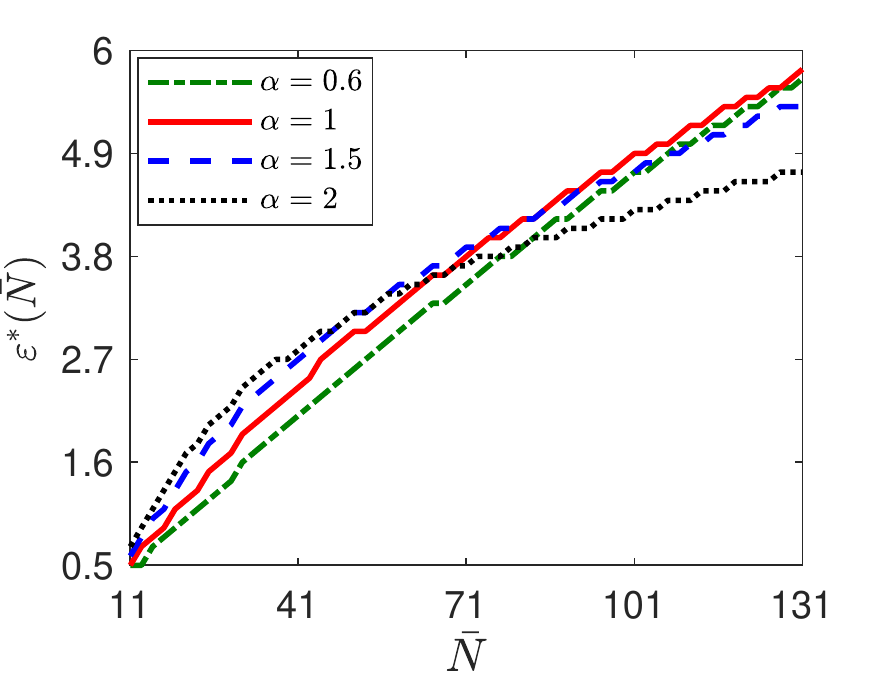}\hspace{-5mm}
\includegraphics[height = 5.76cm, width = 7.86cm]{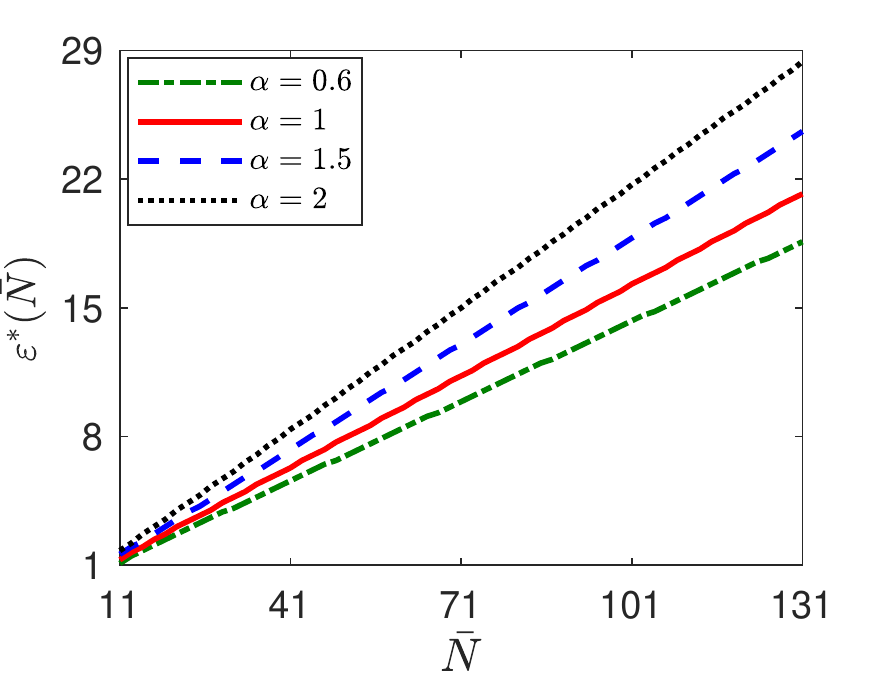}}
\vspace{-1mm}
\centerline{
\includegraphics[height = 5.76cm, width = 7.86cm]{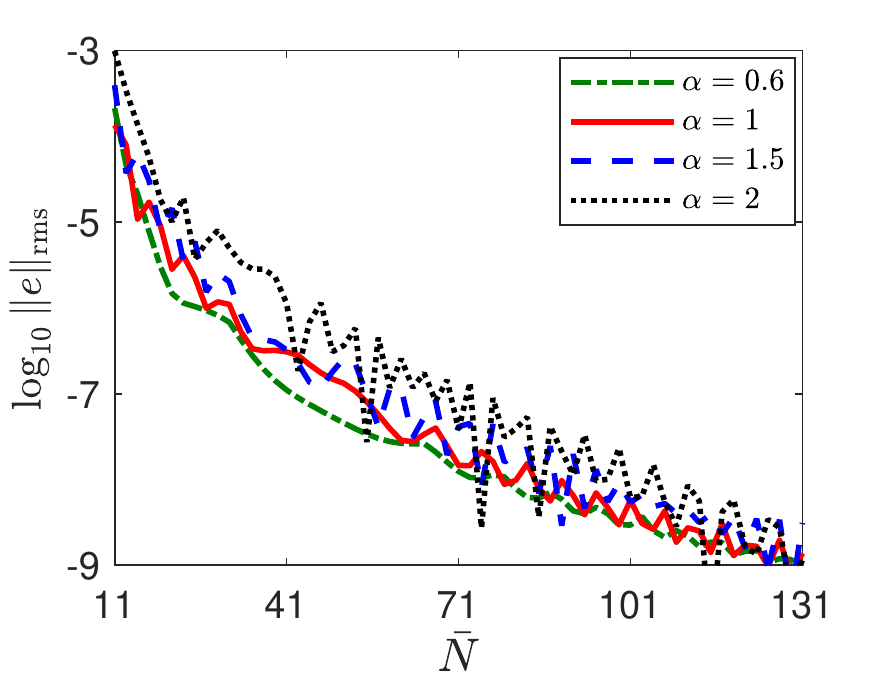}\hspace{-5mm}
\includegraphics[height = 5.76cm, width = 7.86cm]{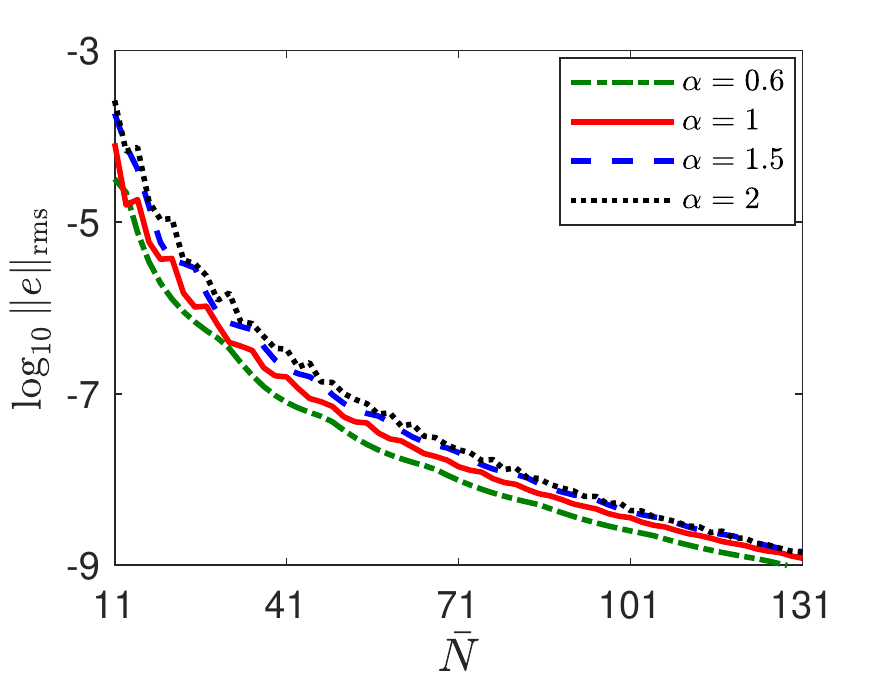}}
\vspace{-1mm}
\centerline{
\includegraphics[height = 5.76cm, width = 7.86cm]{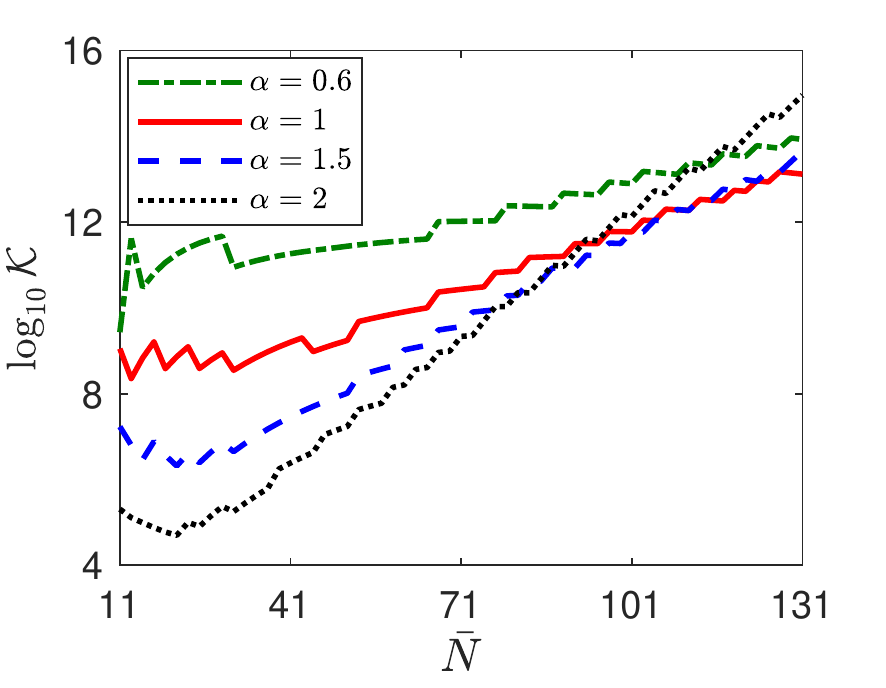}\hspace{-5mm}
\includegraphics[height = 5.76cm, width = 7.86cm]{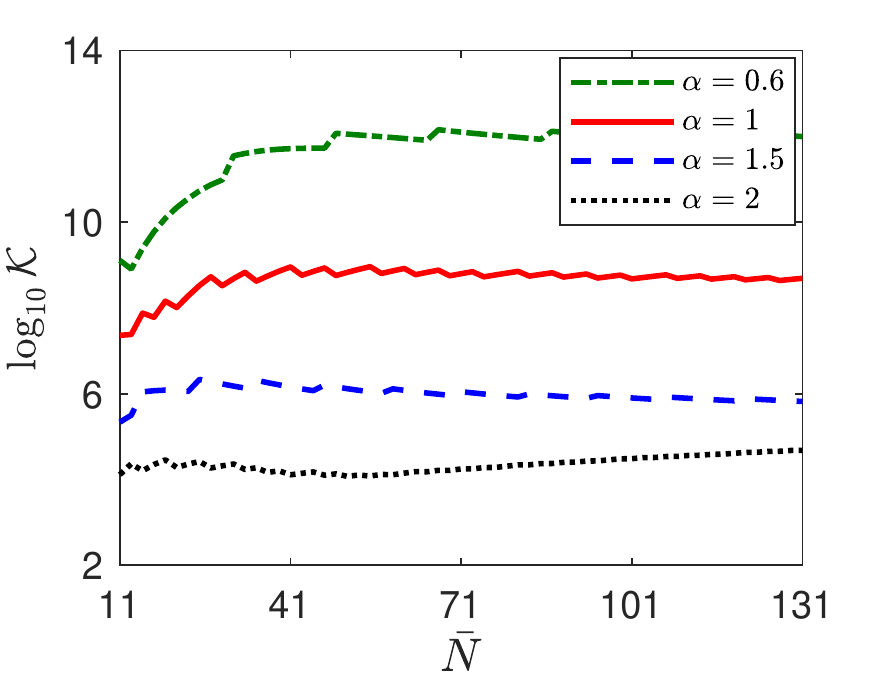}}
\caption{Comparison of the GIMQ-based method (left) and Gaussian-based method in \cite{Burkardt0020} (right) for solving the 1D Poisson problem \eqref{Poisson} with exact solution $u(x) = (1-x^2)^{3+\fl{\ap}{2}}_+$.}\label{Fig4-2-2}
\end{figure} 
In Fig. \ref{Fig4-2-2}, we further compare these two methods by studying their optimal shape parameter $\veps^*(\bar{N})$, and the corresponding numerical errors and condition numbers. 
It shows that for both method, the optimal shape parameter increases almost linearly with  number $\bar{N}$, but  the optimal shape parameters of Gaussian-based method are much larger and sensitive to  power $\ap$. 
More precisely, the larger the power $\ap$, the faster the parameter $\veps^*$ increases with respect to $\bar{N}$. 
The numerical errors at $\veps^*$ are similar for both methods, but the condition numbers of Gaussian-based methods remain almost a constant if $\veps^*$ is used. 
These comparisons could provide some insights for practical applications of these two methods. 

\subsection{Selection of shape parameters}
\label{section4-3}

As we noted, the shape parameter $\veps$ plays an important role in determining the condition number and numerical accuracy of RBF-based methods.  
If a small constant shape parameter (i.e., $\veps_i \equiv \veps$ for $1 \le i \le \bar{N}$) is used, the resulting linear system could easily become ill-conditioned as the number of points $\bar{N}$ increases. 
While using a large shape parameter could suppress the condition number, but at the same time it might deteriorate numerical accuracy if the shape parameter is too big. 
It was pointed out in \cite{Rippa1999} that the optimal shape parameter depends on many factors, including the number of points, distribution of points, properties of functions to approximate, and even computer precision. More discussion can be found in \cite{Tsai2010, Sarra2009, Fornberg2015} and references therein. 
Clearly, it is challenging to find the optimal shape parameters in practice.

In the following, we study two approaches for choosing shape parameters with the goal of balancing numerial accuracy and condition number. 
In the first approach, we consider a constant shape parameter for all RBF center points, which is selected by controlling the condition number in a desired range and thus eliminates the ill-conditioning issues. 
In the second approach, we consider (randomly) different shape parameters for each RBF points. 
Similar approaches can be found in \cite{Sarra2009}, but no studies have been reported in solving fractional PDEs. 

\vskip 8pt
{\bf Condition number indicated shape parameter (denoted as $\veps_{\mathcal K}$). } 
Our studies show that when a constant shape parameter is used,   numerical errors reduce as the number of points increases, and at the same time the condition number quickly increases (e.g. see Tables \ref{Tab4-1-1}--\ref{Tab4-1-2}). 
Consequently, the system could become ill-conditioned if the number of points is too big. 
On the other hand, we find when the number $\bar{N}$ is large, numerical solutions tend to be more accurate if the condition number is large but before the system becomes ill-conditioned. 
Inspired by these observations, we will choose a constant shape parameter such that the resulting condition number falls in a desired range, namely {\it condition number indicated shape parameter}. 
The constant shape parameter selected from this approach will be denoted as $\veps_{\mathcal K}$. 

\begin{figure}[htb!]
\centerline{\includegraphics[height = 5.76cm, width = 7.86cm]{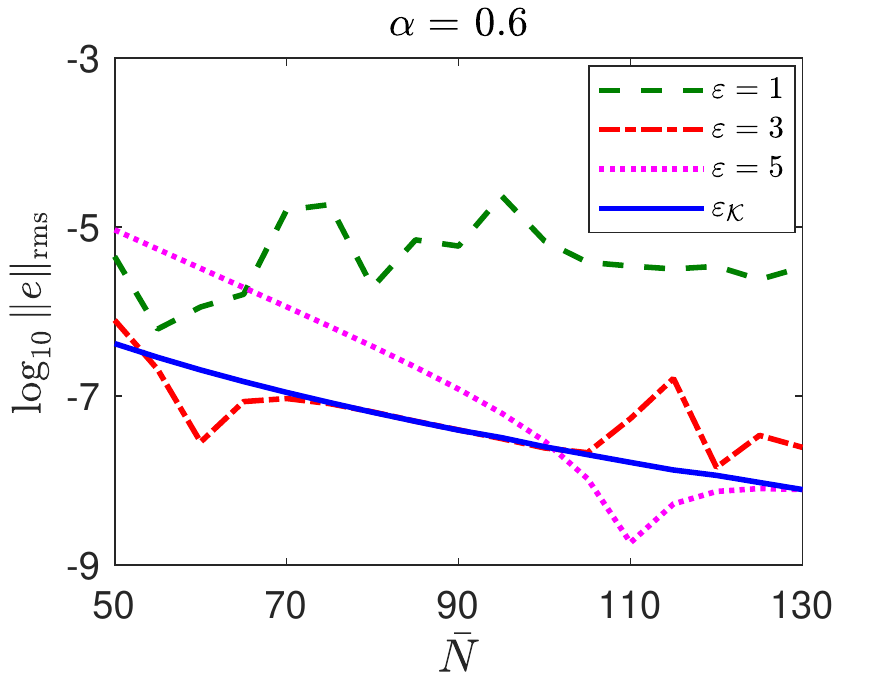}\hspace{-5mm}
\includegraphics[height = 5.76cm, width = 7.86cm]{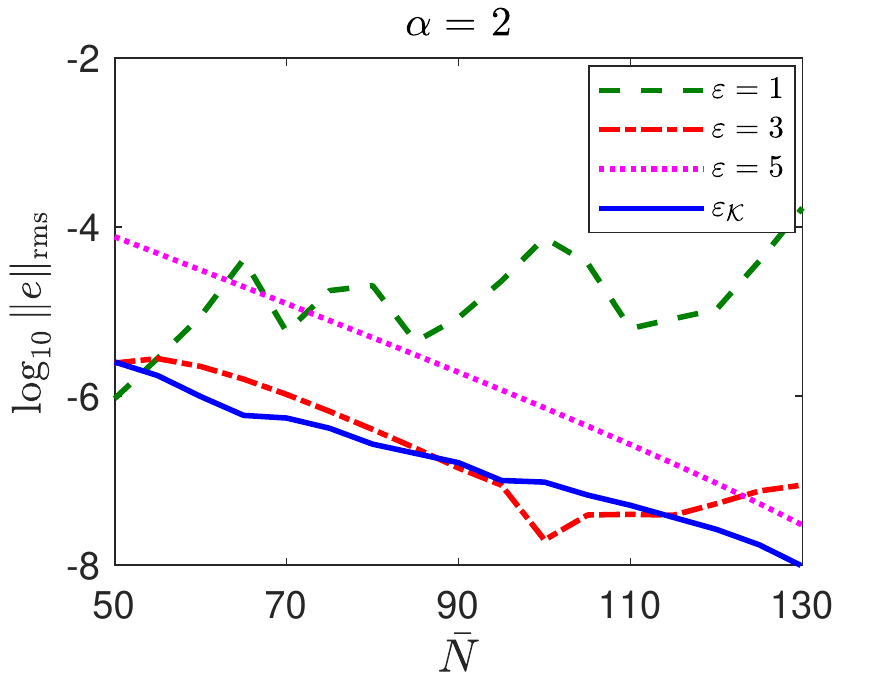}}
\caption{Comparison of numerical errors for different shape parameters in solving the 1D Poisson problem \eqref{Poisson} with exact solution $u(x) = (1-x^2)^{3+\fl{\ap}{2}}_+$.  }\label{Fig4-3-1}
\end{figure}

In Fig. \ref{Fig4-3-1}, we compare numerical errors from different constant shape parameters in solving the Poisson problem \eqref{Poisson} with exact solution $u(x) = (1-x^2)^{3+\fl{\ap}{2}}_+$. 
To determine  $\veps_{\mathcal K}$, we start with a small value of $\veps$ and then increase it gradually until the condition number of the system falls into the range of  $10^{13} \leq \mathcal{K} \leq 10^{16}$,  and the corresponding shape parameter will be defined as  $\veps_{\mathcal K}$. 
Note that the shape parameter selected in this approach changes for different number $\bar{N}$, i.e. $\veps_{\mathcal K} = \veps_{\mathcal K}(\bar{N})$. 
From Fig. \ref{Fig4-3-1}, we find that a constant shape parameter (e.g. $\veps = 1$ or $3$) eventually leads to the growth of numerical errors when the number of points $\bar{N}$ is large enough because the system becomes ill-conditioned. 
Moreover, a large shape parameter always lead to bigger numerical errors even though the system is well-conditioned (see, e.g., $\veps = 5$).  
In contrast,  the ${\mathcal K}$-indicated shape parameter $\veps_{\mathcal K}$ can well balance the condition number and numerical accuracy.  
It can prevent the system from being ill-conditioned and ensure small numerical errors as $\bar{N}$ increases. 
This property is important for practical simulations. 


\vskip 10pt
{\bf Random-perturbed variable shape parameter (denoted as $\veps_\dt$). }  
The condition number indicated shape parameter $\veps_{\mathcal K}$ can effectively avoid the ill-conditioning issues in simulations with large $\bar{N}$. 
Next, we will study another approach --  {\it random-perturbed variable shape parameter}. 
It perturbs a constant shape parameter with uniformly distributed random numbers at RBF center points, and thus  the shape parameter can be viewed as a random function of center point $\bx_i$ \cite{Sarra2009b, Wertz2006}. 
Let $[\veps_{\min},  \veps_{\max}]$ denote the interval that the shape parameter is allowed. 
The random-perturbed shape parameter $\veps_{\dt, i}$ for center point $\bx_i$ is given by:
\bea\label{fun3-3-1} 
\varepsilon_{\dt, i} = \varepsilon_{\rm min} + \dt_i (\varepsilon_{\rm max} - \varepsilon_{\rm min}),
\eea
where $\dt_i$ is a uniformly distributed random number on the interval (0,1).  

\begin{figure}[htb!]
\centerline{\includegraphics[height = 5.76cm, width = 7.86cm]{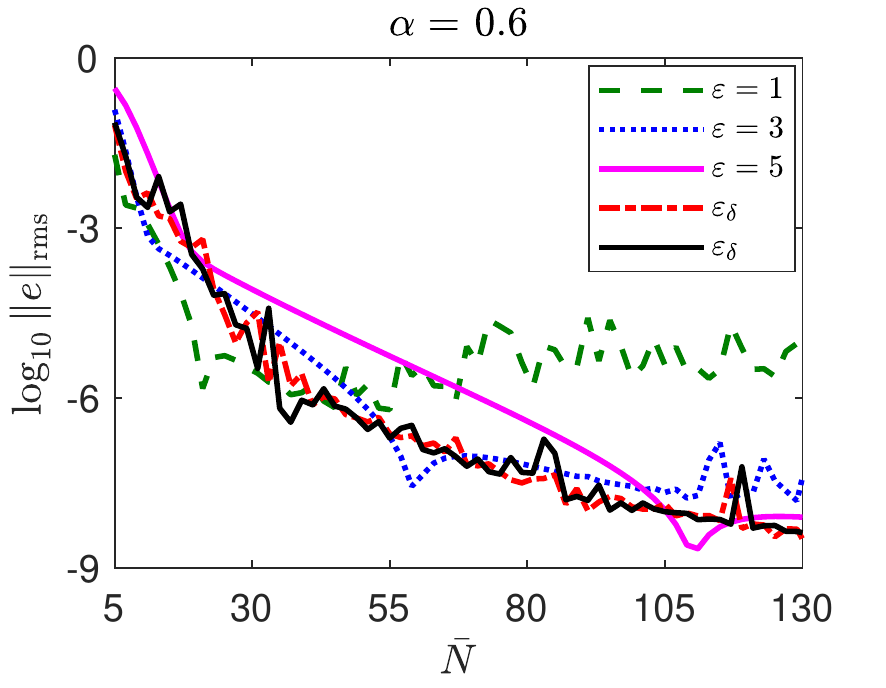}\hspace{-5mm}
\includegraphics[height = 5.76cm, width = 7.86cm]{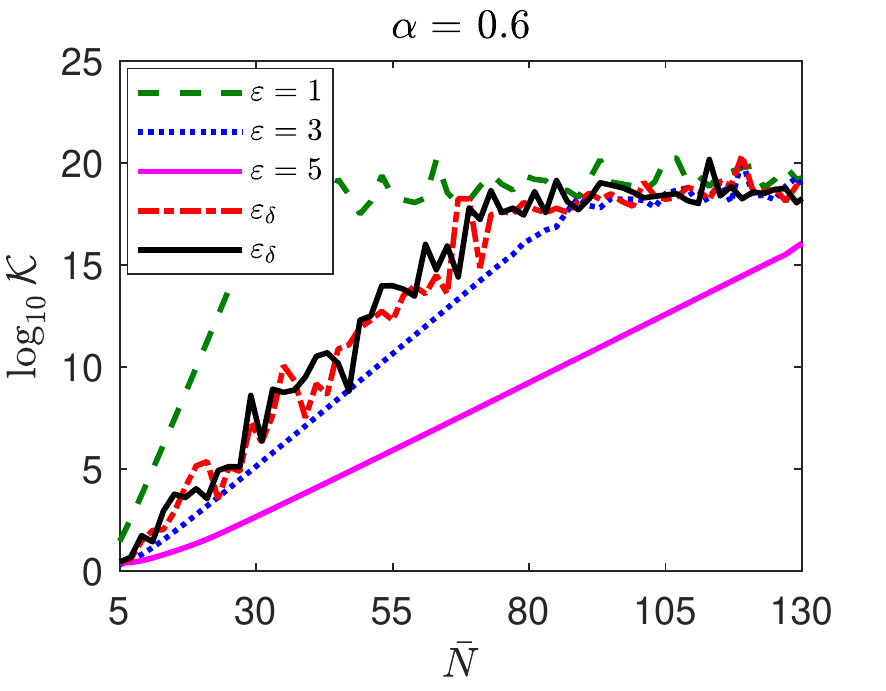}}
\centerline{\includegraphics[height = 5.76cm, width = 7.86cm]{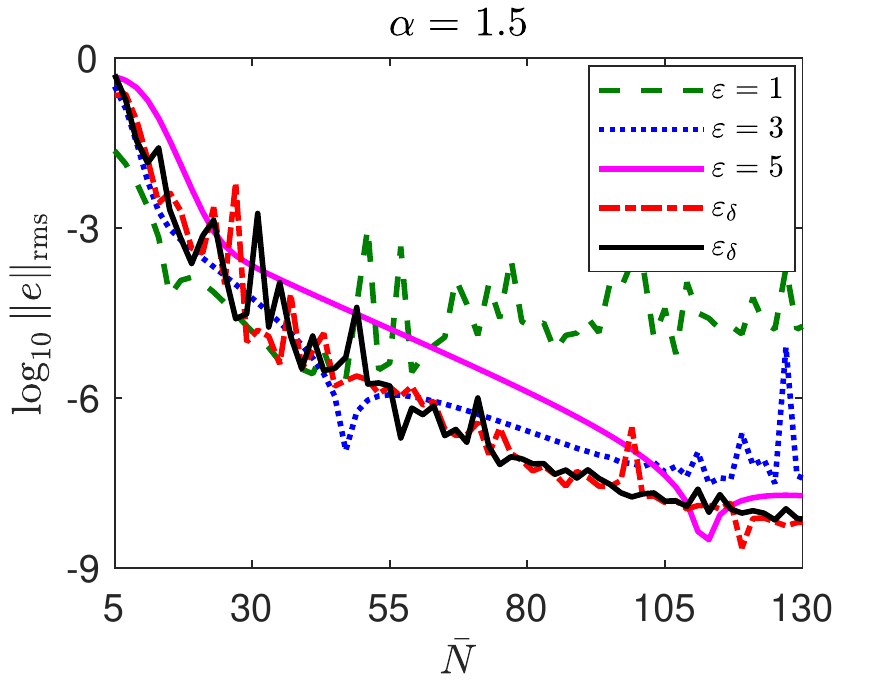}\hspace{-5mm}
\includegraphics[height = 5.76cm, width = 7.86cm]{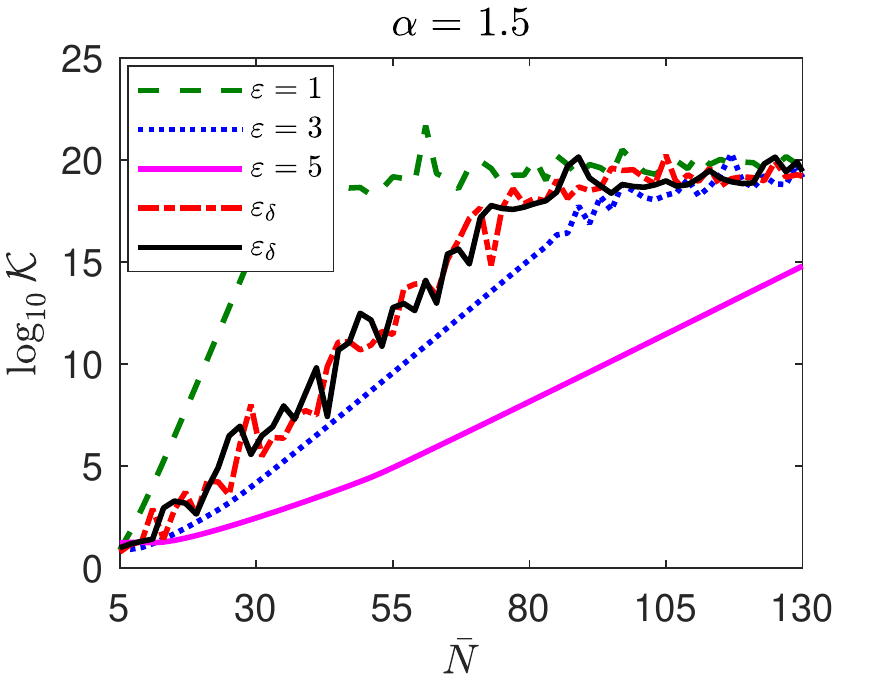}}
\caption{Comparison of numerical errors and condition numbers for different shape parameters in solving the 1D Poisson problem \eqref{Poisson} with exact solution $u(x) = (1-x^2)^{3+\fl{\ap}{2}}$.}\label{Fig4-3-2}
\end{figure}
In Fig. \ref{Fig4-3-1}, we test the performance of random-perturbed variable shape parameters $\veps_\dt$ in comparison to constant parameters $\veps$, where the Poisson problem (\ref{Poisson}) with exact solution $u(x) = (1-x^2)^{3+\fl{\ap}{2}}_+$ is studied. 
We choose $\varepsilon_{\min}=1$ and $\varepsilon_{\max}=5$ in \eqref{fun3-3-1}. 
It shows that two simulations with different $\veps_\dt$ give similar numerical errors and condition numbers. 
Compared to constant shape parameters, using distinct shape parameters for each RBF center point helps suppressing the condition number and  leads to good numerical errors for different number $\bar{N}$. 
It can also save the computational time in searching for a proper shape parameter  (e.g., $\veps_{\mathcal K}$).  
Here, the choices of [$\veps_{\min}$, $\veps_{\max}$] is important in determining random-perturbed variable shape parameters. 
Our studies in Section \ref{section4-2} suggest that this range might depend on the minimum distance between RBF center points, which we will further explore in our future study.

\section{More applications}
\setcounter{equation}{0}
\label{section5}

In this section, we will focus on our GIMQ-based meshless method and apply it to study both elliptic problems and diffusion equations.  
In our simulations, we will use the random-perturbed shape parameters as described in \eqref{fun3-3-1}, and  the test points are chosen from the same set of RBF center points (consequently, $M = N$ in \eqref{scheme5}). 

\subsection{Elliptic problem on an L-shaped domain}
\label{section5-1}

We consider the following elliptic problem on an L-shaped  domain, i.e., $\Og = \{(x, y) \in (-1, 1)^2\backslash [0, 1)^2\}$, with nonhomogeneous Dirichlet boundary conditions: 
\bea\label{fun4-1-1}
(-\Dt)^{\fl{\ap}{2}} u + 2u = f, \quad \mbox{for}  \ \ \bx\in \Omega; \qquad\ u(\bx) = e^{-|\bx|^2}, \quad  \mbox{for} \ \  \bx\in \Upsilon. 
\eea
The right hand side  function $f$ in (\ref{fun4-1-1}) is chosen as:
\bea\label{fun4-1-2}
f(\bx) =  2^\ap \Gamma(1+\fl{\ap}{2})\,_1F_1\Big(1+\fl{\ap}{2}; \, 1; \, -|\bx|^2\Big) + 2e^{-|\bx|^2}
\eea
where $_1F_1$ denotes the confluent hypergeometric function, and the exact solution  is  $u = e^{-|\bx|^2}$ for $\bx \in {\mathbb R}^2$. 
We choose RBF center/test points as uniformly distributed grid points on domain $\bar{\Og}$, and $\veps_{\rm min} = 0.1$ and $\veps_{\rm max} = 4$ in the random-perturbed shape parameters (\ref{fun3-3-1}). 

Table \ref{Tab5-1-1} shows numerical errors and condition numbers for various $\ap$.  
It is evident that numerical errors decrease with a spectral rate as $\bar{N}$ increases. 
\begin{table}[htb!]
\begin{center} 
\begin{tabular}{|c|c|c|c|c|c|c|c|c|}
\hline
&\multicolumn{2}{|c|}{$\ap = 0.6$} & \multicolumn{2}{|c|}{$\ap = 1$} & \multicolumn{2}{|c|}{$\ap = 1.5$} & \multicolumn{2}{|c|}{$\ap = 2$}   \\
\cline{2-9}
$\bar{N} $ & RMS & $\mathcal{K}$ & RMS & $\mathcal{K}$ & RMS & $\mathcal{K}$ & RMS & $\mathcal{K}$ 
 \\ 
\hline
$21$ & 6.238e-2 & 8.909e2 & 4.888e-2 & 1.147e2 & 3.607e-2 & 1.288e2 & 2.339e-2 & 2.036e2   \\
\hline
$65$ & 3.795e-3 & 2.227e7 & 8.403e-3 & 7.030e6 & 4.834e-3 & 2.901e6 & 9.294e-3 & 1.794e5  \\
\hline
$133$ & 2.892e-4 & 2.420e10 & 3.976e-4 & 1.234e10 & 1.779e-4 & 9.662e9 & 1.702e-5 & 1.199e10  \\
\hline
$225$ & 3.899e-6 & 4.069e11 & 1.367e-6 & 7.021e11 & 4.178e-5 & 8.393e14 & 2.626e-6 & 1.236e11  \\
\hline
$341$ & 3.902e-8 & 1.845e15 & 1.336e-8 & 3.219e15 & 4.239e-8 & 3.178e17 & 2.159e-8 & 1.481e17 \\
\hline
\end{tabular}
\caption{Numerical errors  $\|e\|_{\rm rms}$ and condition number $\mathcal{K}$ in solving the elliptic problem (\ref{fun4-1-1}) on an L-shaped domain. }\label{Tab5-1-1}
\end{center}
\end{table}
Compared to finite difference/element methods, our method can yield a more accurate result with fewer number of points. 
For example, to achieve an accuracy of ${\mathcal O}(10^{-8})$, our method only requires the number of points $\bar{N} = 341$ (equivalently, distance between two grid points $h = 0.1$). 
Due to the nonlocality of the fractional Laplacian, the discretization of fractional PDEs usually results in a linear system with full matrix; consequently their simulations require huge storage and computational costs. 
While our method yielding higher accuracy with fewer points can effectively relax these challenges in  simulations of fractional PDEs. 


\subsection{Normal and anomalous diffusion problems}
\label{section5-2}

We solve the diffusion problem with $f \equiv 0$ and $\kappa = 0.5$ in \eqref{diffusion},  and homogeneous Dirichlet boundary conditions in \eqref{diffusion-BC}. 
Set the domain  $\Og = (-2, 2)^2\backslash[0.5, 1.5]^2$, and the initial condition 
\beas
u(\bx, 0) = 3e^{-10|\bx|^4}, \qquad \mbox{for} \ \ \bx \in {\mathbb R}^2.
\eeas
We then compare normal ($\ap = 2$) and anomalous ($\ap < 2$) diffusion by studying time evolution of the  solution $u(\bx, t)$. 
In our simulations, RBF center/test points are chosen as uniformly distributed grid points on $\bar{\Og}$. 
The random-perturbed shape parameters are used with  $\veps_{\rm min} = 1$ and $\veps_{\rm max} = 5$.

\begin{figure}[htb!]
\centerline{\includegraphics[height = 4.6cm, width = 5.86cm]{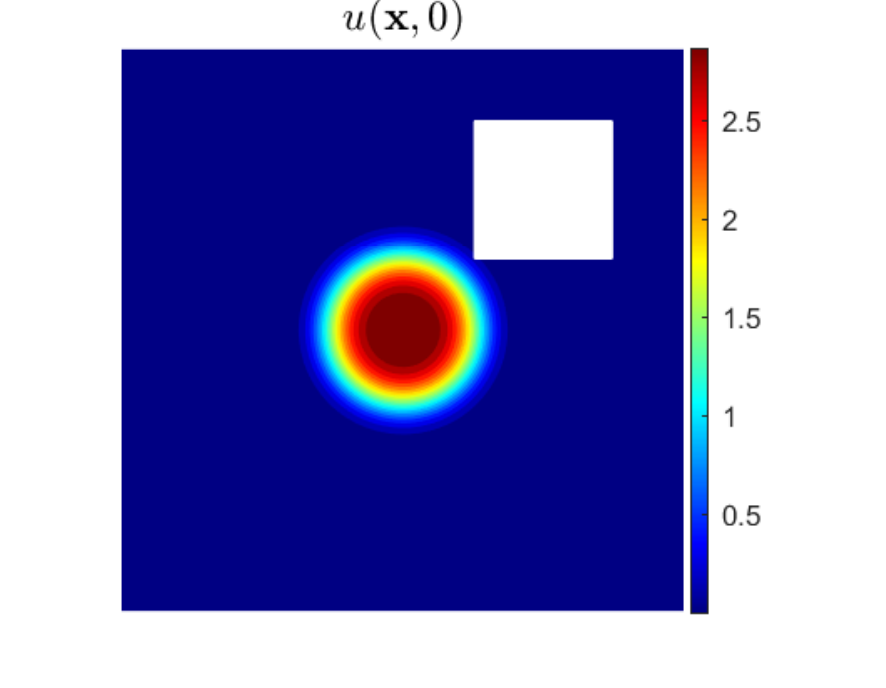}
\includegraphics[height = 5.26cm, width = 7.36cm]{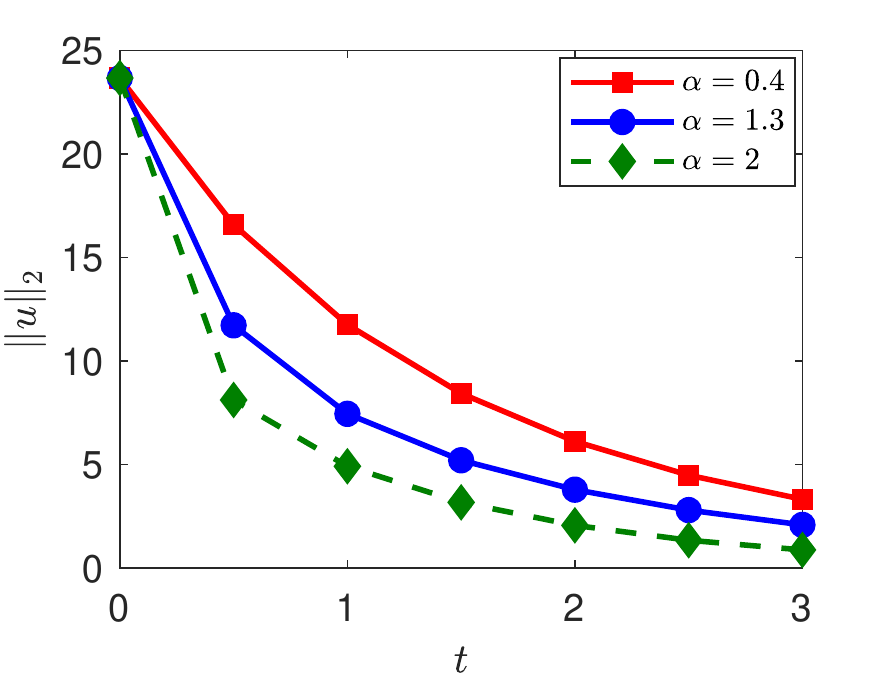}}
\caption{(a) Initial condition $u(\bx, 0)$ of the diffusion problem; (b) Dynamics of $\|u(\cdot, t)\|_2$ in normal and anomalous diffusion.}\label{Fig5-2-1}
\end{figure}

Fig. \ref{Fig5-2-1} shows the initial solution $u(\bx, 0)$ and dynamics of  solution norm $$\|u(\cdot, t)\|_2 = \Big(\int_{\Og}|u(\bx, t)|^2d\bx\Big)^{\fl{1}{2}}, \qquad t \ge 0$$ for different power $\ap$. 
Here, we choose time step $\tau = 0.001$, and the number of points ${\bar{N}} = 600$. 
We have verified and ensured the convergence of solutions with smaller time step and more RBF points. 
It shows that initially the radially-symmetric solution concentrates around $\bx = {\bf 0}$ with norm 
$\|u(\cdot, 0)\|_2 = 23.67$. 
Then norm $\|u(\cdot, t)\|_2$ monotonically decreases over time owing to zero source term (i.e., $f = 0$) and homogeneous Dirichlet boundary conditions.  
It shows that the smaller the power $\ap$, the slower the decay of solution norm $\|u(\cdot, t)\|_2$.  

\begin{figure}[htb!]
\centerline{\includegraphics[height = 3.86cm, width = 5.16cm]{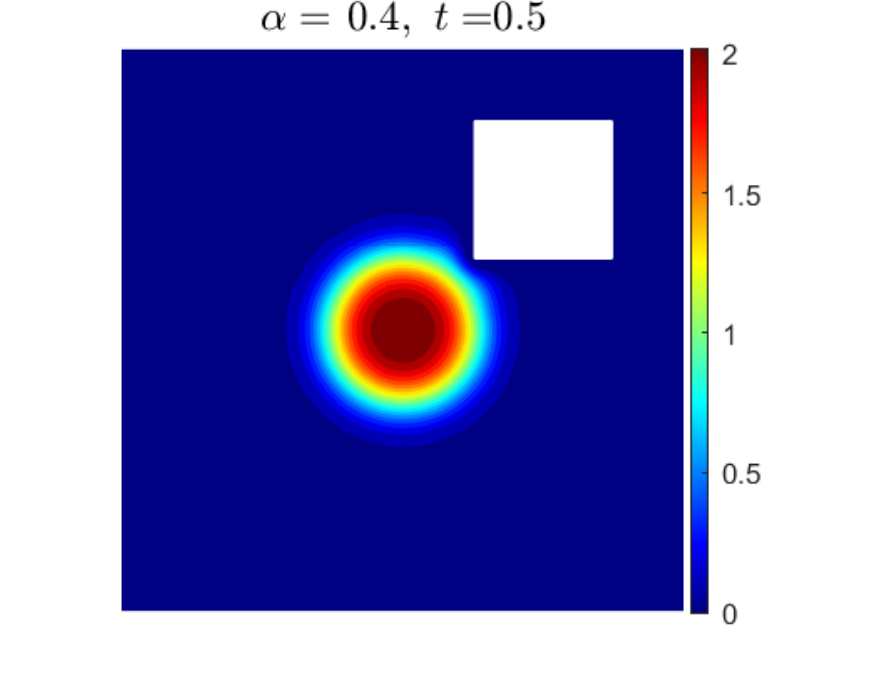}\hspace{-5mm}
\includegraphics[height = 3.86cm, width = 5.16cm]{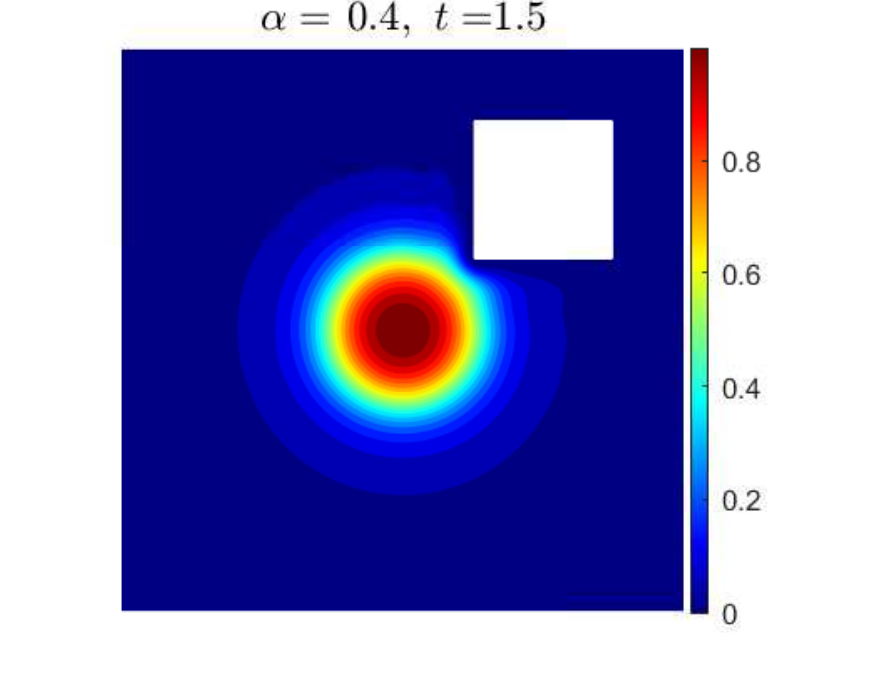}\hspace{-5mm}
\includegraphics[height = 3.86cm, width = 5.16cm]{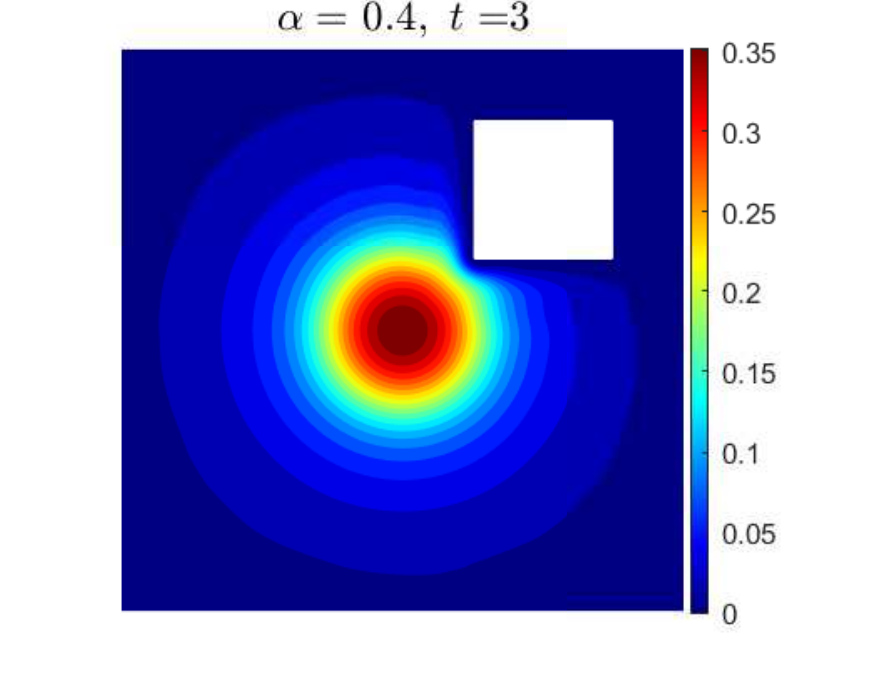}}
\centerline{
\includegraphics[height = 3.86cm, width = 5.16cm]{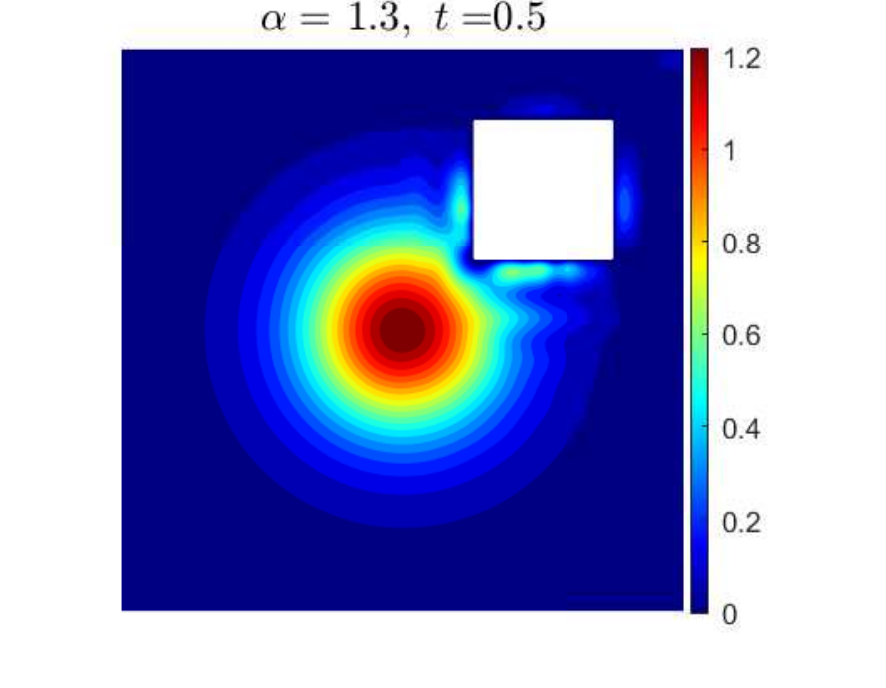}\hspace{-5mm}
\includegraphics[height = 3.86cm, width = 5.16cm]{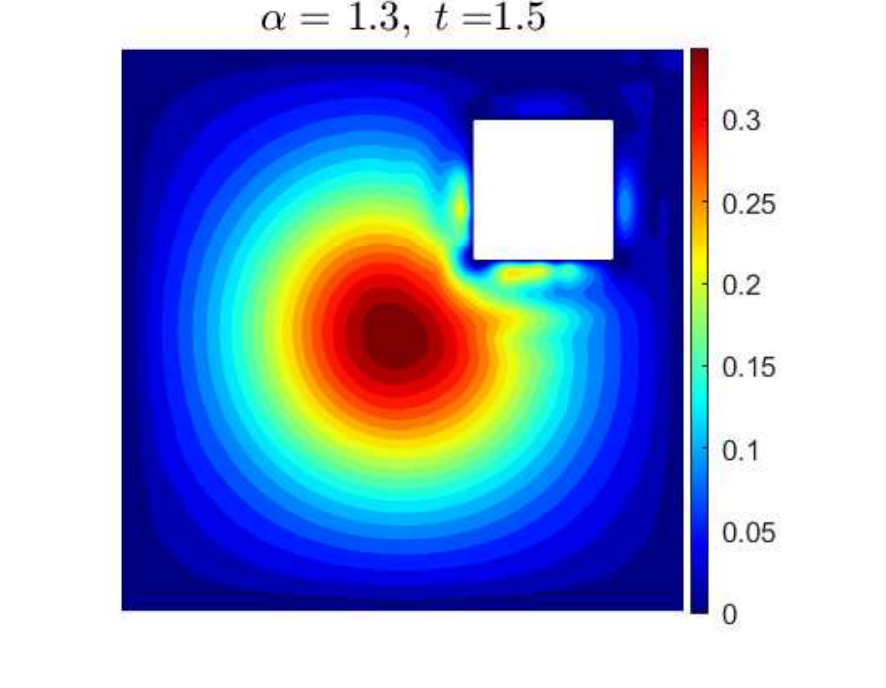}\hspace{-5mm}
\includegraphics[height = 3.86cm, width = 5.16cm]{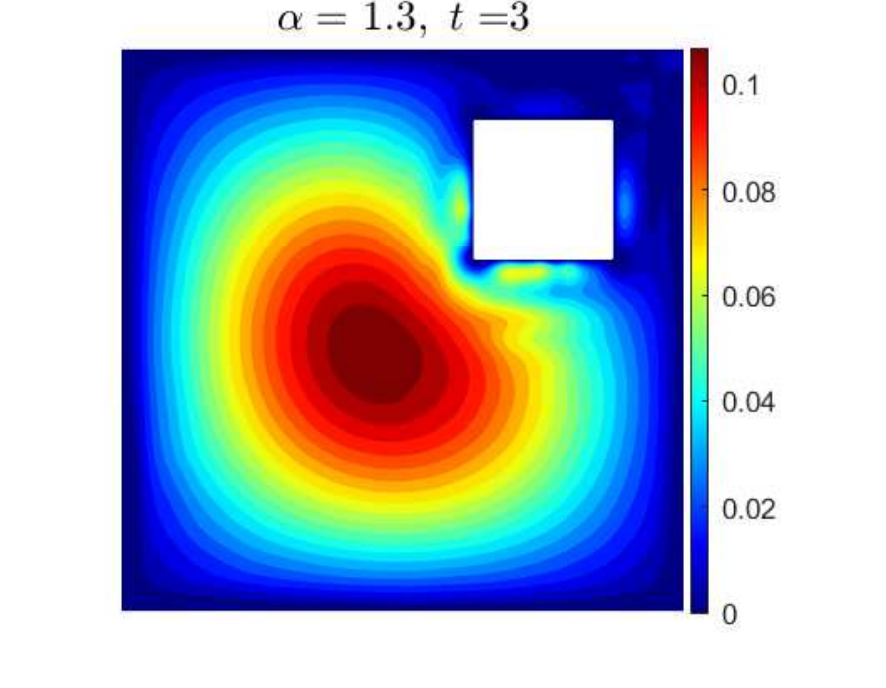}}
\centerline{
\includegraphics[height = 3.86cm, width = 5.16cm]{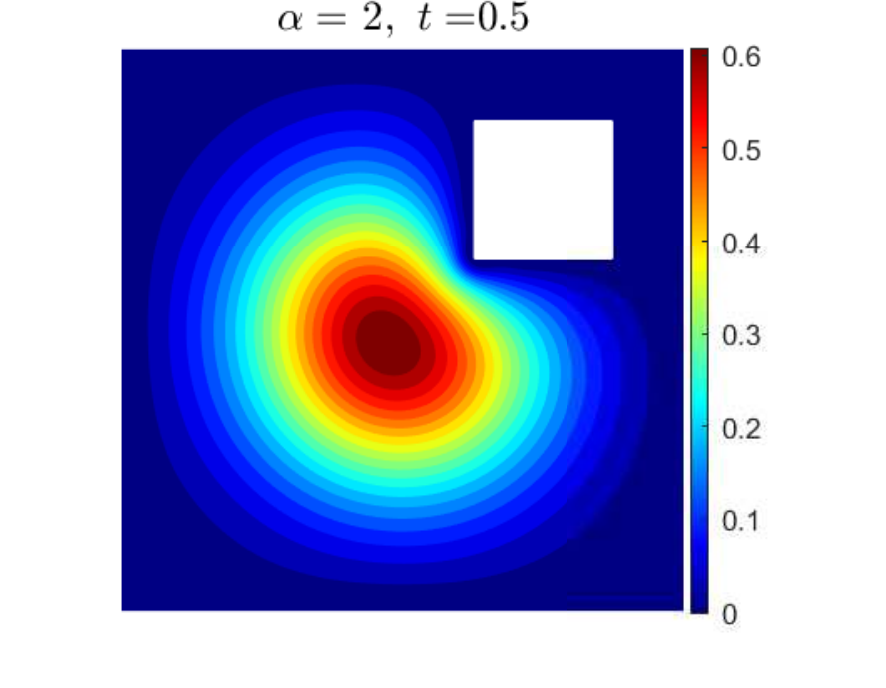}\hspace{-5mm}
\includegraphics[height = 3.86cm, width = 5.16cm]{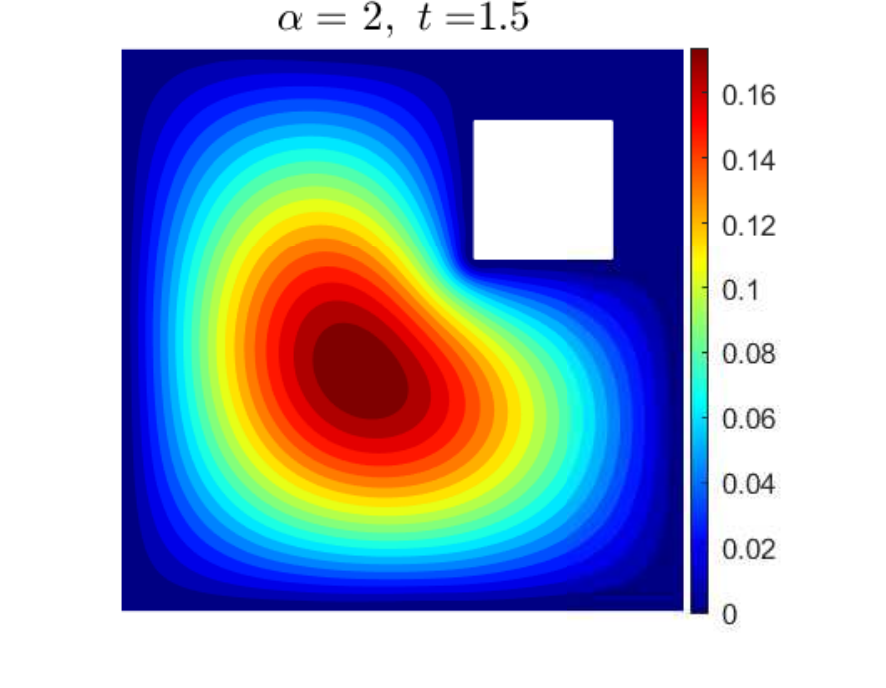}\hspace{-5mm}
\includegraphics[height = 3.86cm, width = 5.16cm]{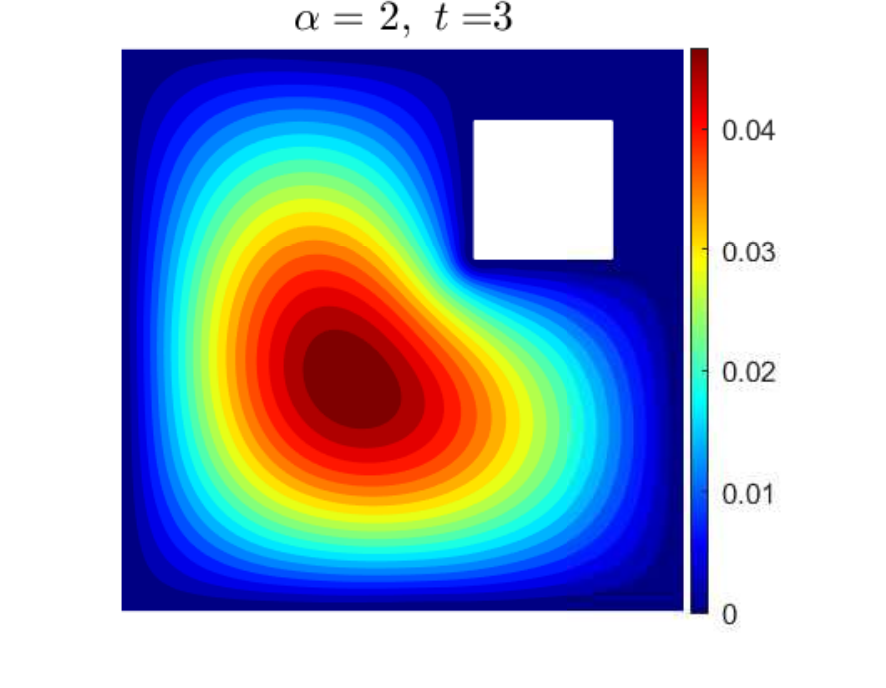}}
\caption{Time evolution of solution in normal and anomalous diffusion problems.}
\label{Fig5-2-2}
\end{figure}

Fig. \ref{Fig5-2-2} further compares the solution evolution at different time $t$. 
It shows that  solution  $u(\bx, t)$ diffuses  from domain center towards boundaries over time, and it starts deforming once reaching the inner boundaries (i.e., boundary of region $[0.5, 1.5]^2$). 
It is clear that the smaller the power $\ap$, the slower the solution diffuses, consistent with the slower decay of norm $\|u(\cdot, t)\|_2$. 
At time $t = 1.5$,  the solution of $\ap = 1.3$ and $2$ have reached the outer boundaries (i.e., $x = \pm 2$ or $y = \pm 2$),  and then the amplitude of solution continues decreasing due to the zero boundary conditions. 
In contrast, solution of $\ap = 0.4$ reaches the inner boundary at time $t = 1.5$, but is still far away from the outer boundaries until $t = 3$. 

\subsection{Heat conduction problem}
\label{section5-3}

In this example, we continue our study on normal and anomalous diffusion and explore the boundary effects in classical and fractional PDEs. 
Consider the problem 
\bea\label{fun4-2-1}
\p_t u(\bx, t) = -\kappa(-\Dt)^\fl{\ap}{2}u + u, &\ \, &\mbox{for} \ \ \bx\in\Og, \quad t > 0, 
\eea 
where the domain  $\Og = (-1, 1)^2$. 
Let the initial condition $u(\bx, 0) \equiv 0$ for $\bx \in \bar{\Og}$. 
The boundary conditions are set as 
\beas
u(\bx, t) = \left\{\begin{array}{ll}
\sin\big[\pi(x-x_c+\fl{1}{2})\big]\sin\big[\fl{\pi}{2}(y+1)\big], \quad \ & \mbox{for} \  \ \bx \in \Xi, \\
0, & \mbox{for} \  \ \bx \in \Upsilon\backslash\Xi,
\end{array}\right. \qquad t \ge 0,
\eeas
where the stripe region $\Xi = [x_c, x_c+\fl{1}{4}]\times[-1,1]$ for $x_c \ge 1$. 
In other words,  boundary conditions are zero everywhere except on  region $\Xi$, and the location of region $\Xi$ is controlled by the value of $x_c$. 
When $x_c = 1$, region $\Xi$ connects to domain $\Og$, while if $x_c > 1$ there is no contact between $\Xi$ and $\Og$. 
This allows use to study the effects of boundary conditions by adjusting the value of $x_c$.

\begin{figure}[htb!]
\centerline{\includegraphics[height = 3.96cm, width = 4.36cm]{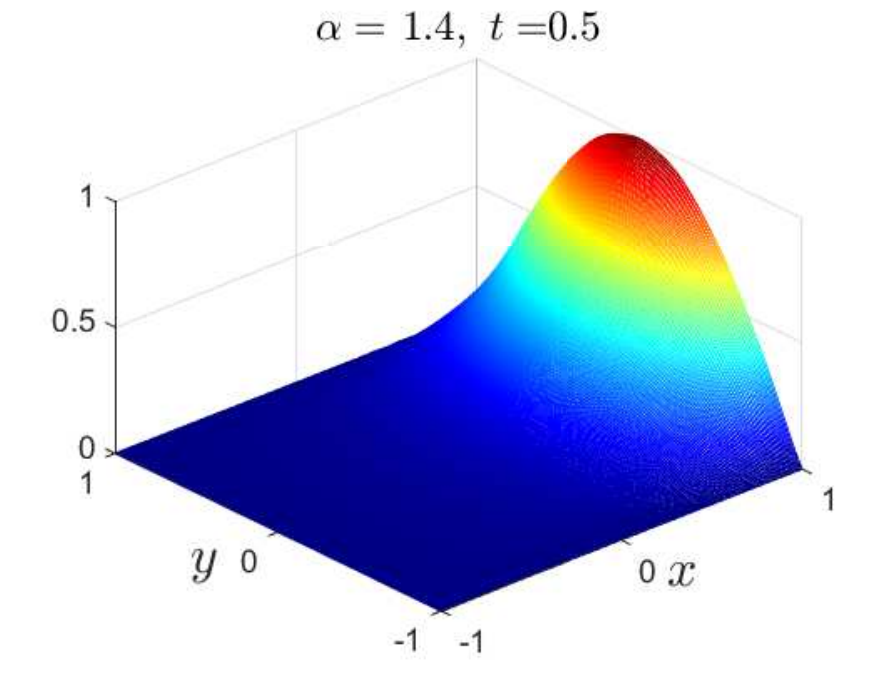}\hspace{-5mm}
\includegraphics[height = 3.96cm, width = 4.36cm]{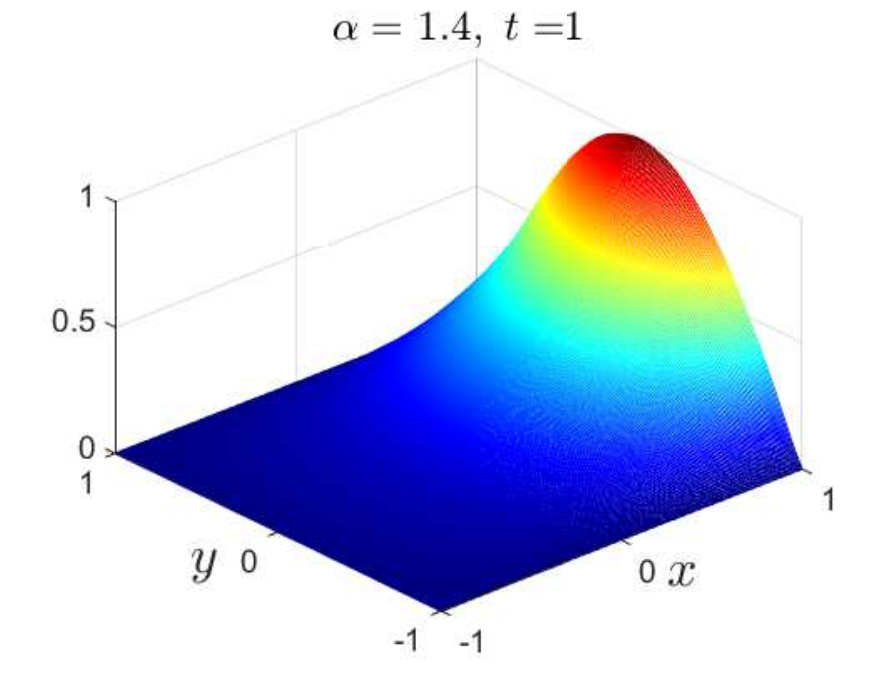}\hspace{-5mm}
\includegraphics[height = 3.96cm, width = 4.36cm]{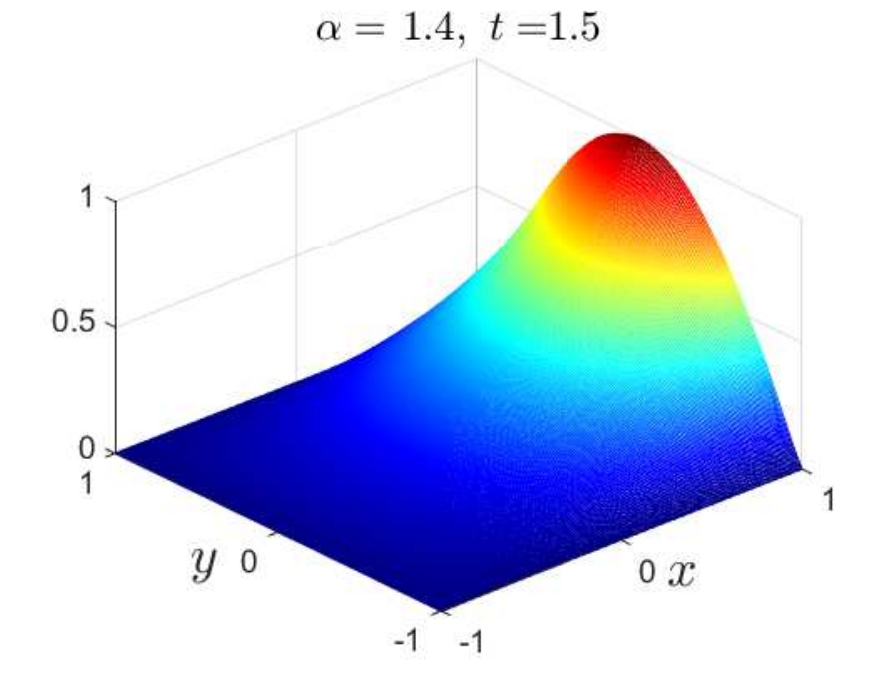}\hspace{-5mm}
\includegraphics[height = 3.96cm, width = 4.36cm]{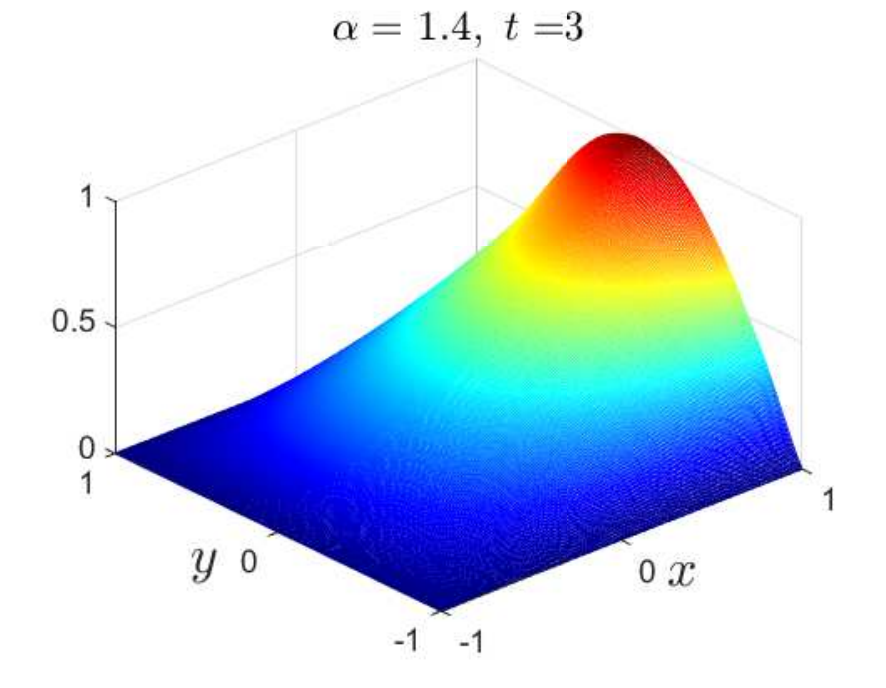}}
\centerline{\includegraphics[height = 3.96cm, width = 4.36cm]{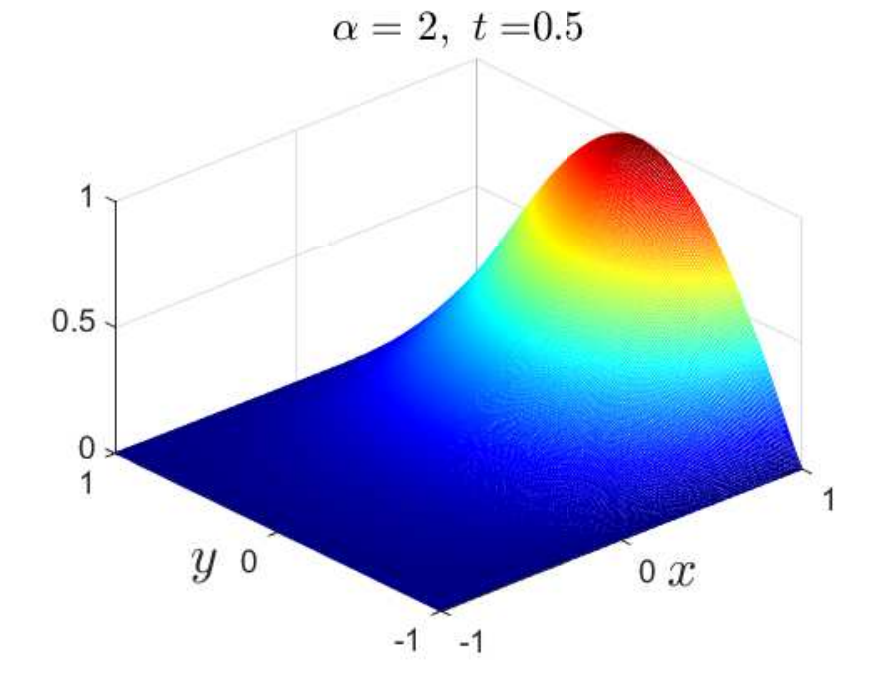}\hspace{-5mm}
\includegraphics[height = 3.96cm, width = 4.36cm]{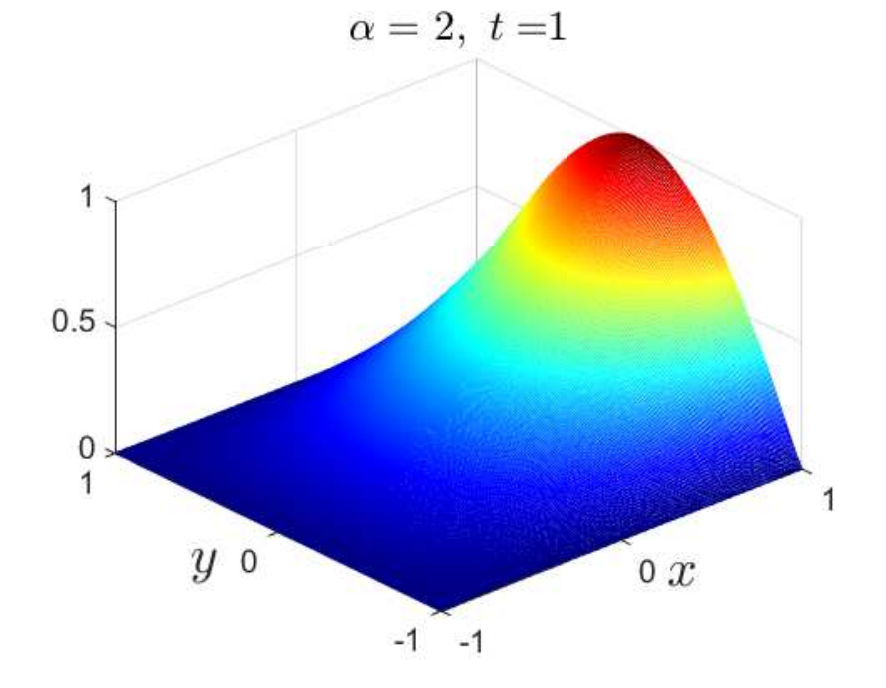}\hspace{-5mm}
\includegraphics[height = 3.96cm, width = 4.36cm]{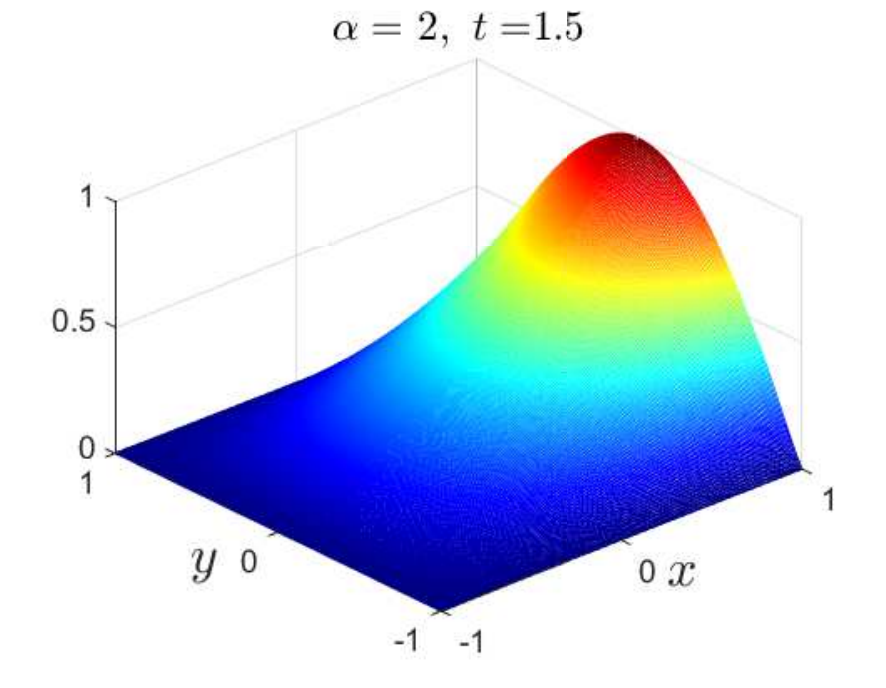}\hspace{-5mm}
\includegraphics[height = 3.96cm, width = 4.36cm]{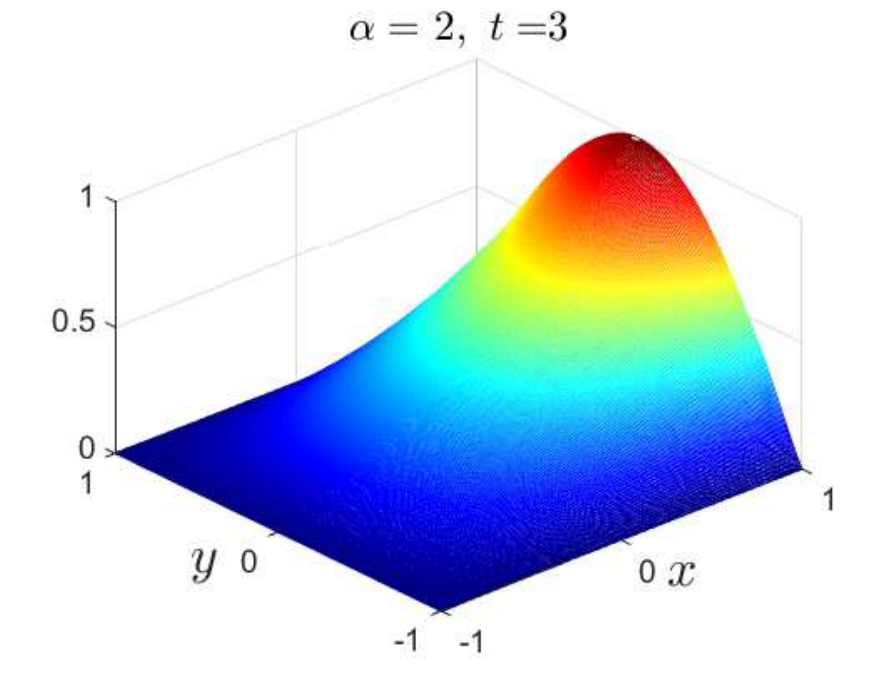}}
\caption{Numerical solution of heat equation \eqref{fun4-2-1} at different time, where $x_c = 1$.  }\label{Fig5-3-1}
\end{figure}

\begin{figure}[htb!]
\centerline{\includegraphics[height = 5.76cm, width = 7.86cm]{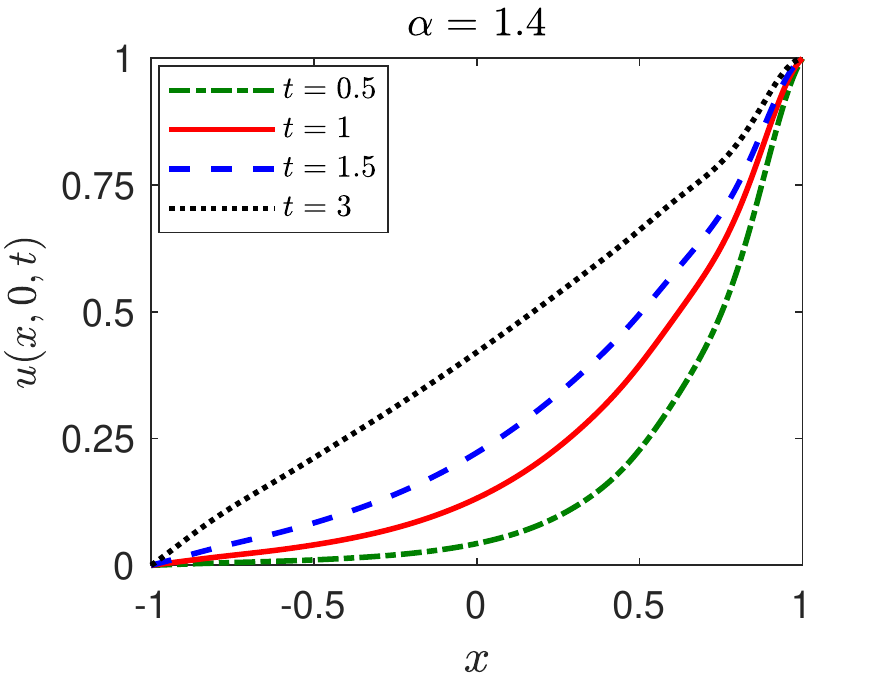}\hspace{-5mm}
\includegraphics[height = 5.76cm, width = 7.86cm]{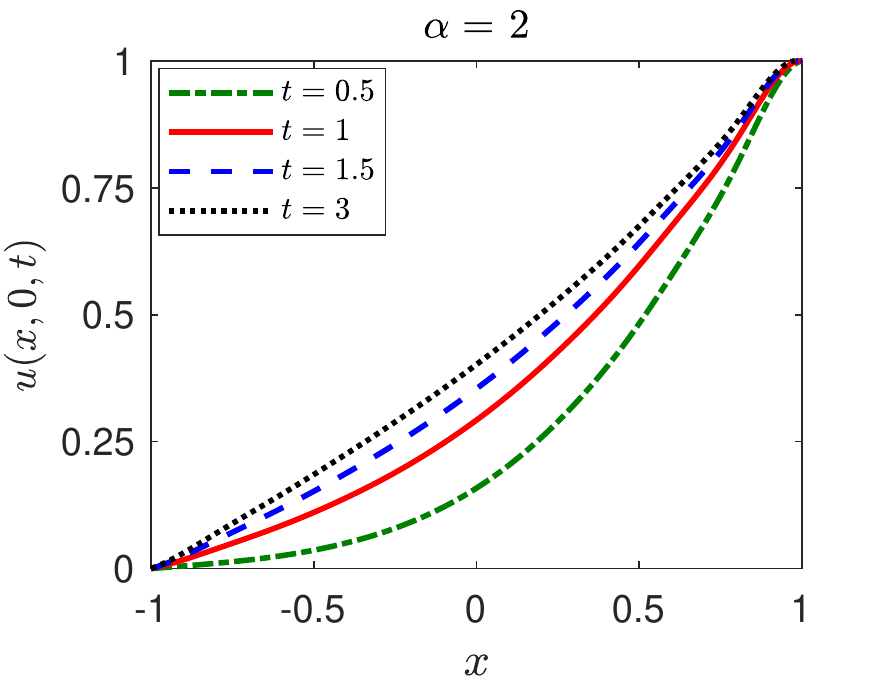}}
\caption{Time snapshot of $u(x, 0, t)$ for $\ap = 1.4$ (a) and $\ap = 2$ (b).}\label{Fig5-3-2}
\end{figure}

Fig. \ref{Fig5-3-1} shows the dynamics of solution $u(\bx, t)$ for $x_c = 1$, and $\ap = 1.4,\, 2$. 
Initially, the solution $u(\bx, 0) = 0$. 
Due to nonzero boundary conditions on  $\Xi = [1, 1.25]\times[-1, 1]$, the solution around $x = 1$ quickly increases and simultaneously diffuses into domain $\Og$. 
It remains symmetric with respect to $y = 0$ for any time $ > 0$, consistent with the symmetry of boundary conditions. 
Note that when $\ap = 2$, the diffusion operator $-\Dt$ is local, and the effective nonzero boundary conditions occur only at line $x = 1$. 
Fig. \ref{Fig5-3-2} further compares the solution evolution at $y = 0$ for different time $t$. 
It shows that the solution of classical cases diffuse fast at the initial stage, and it quickly approaches the steady state. 

\begin{figure}[htb!]
\centerline{
\includegraphics[height = 3.96cm, width = 4.36cm]{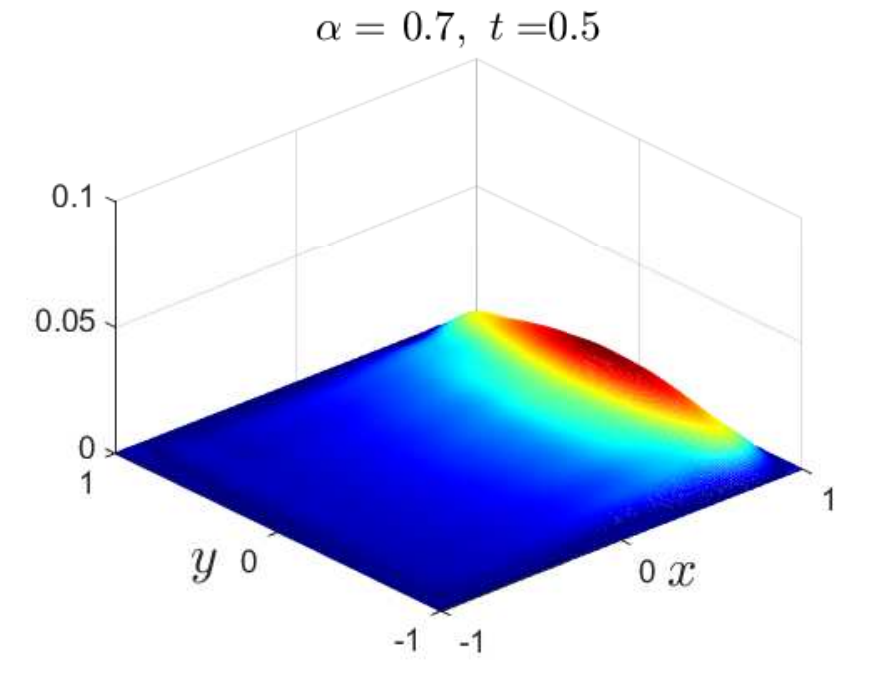}\hspace{-5mm}
\includegraphics[height = 3.96cm, width = 4.36cm]{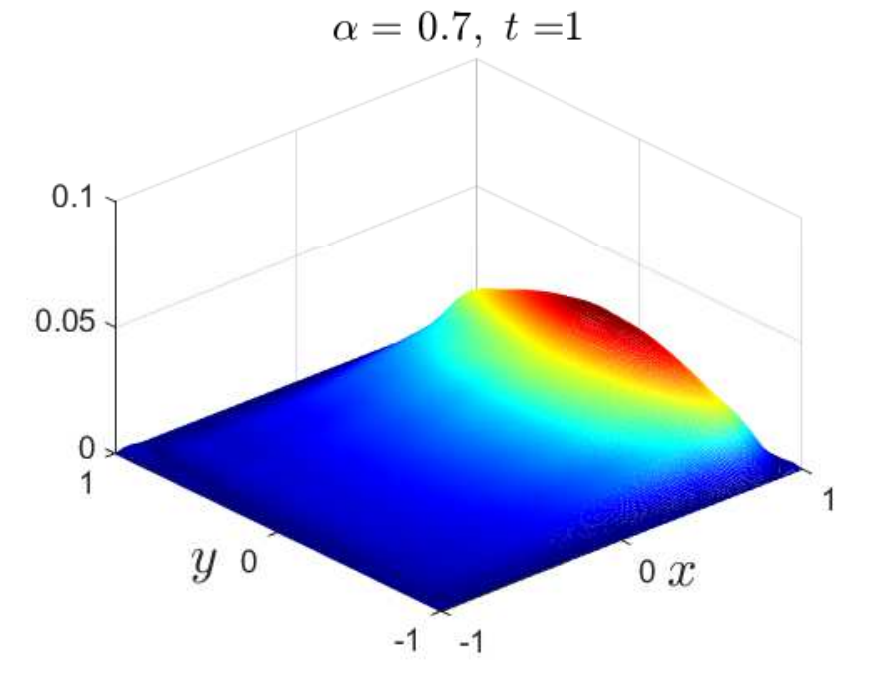}\hspace{-5mm}
\includegraphics[height = 3.96cm, width = 4.36cm]{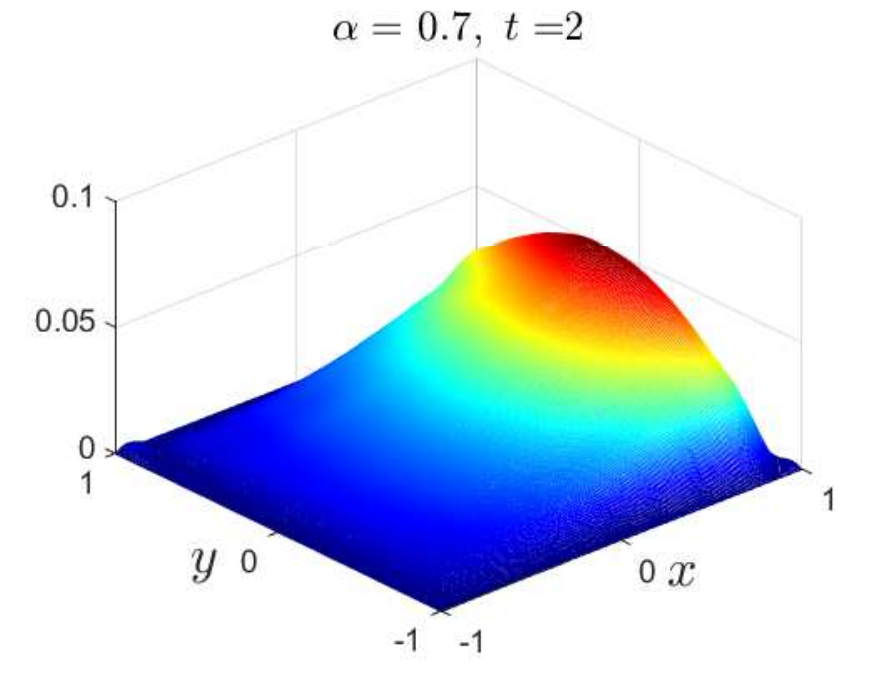}\hspace{-5mm}
\includegraphics[height = 3.96cm, width = 4.36cm]{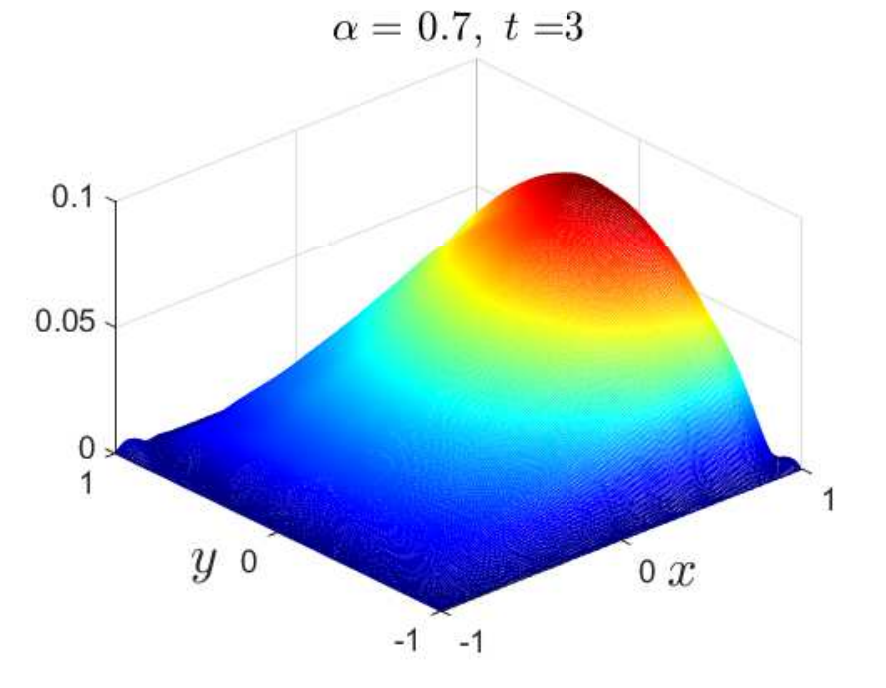}}
\centerline{
\includegraphics[height = 3.96cm, width = 4.36cm]{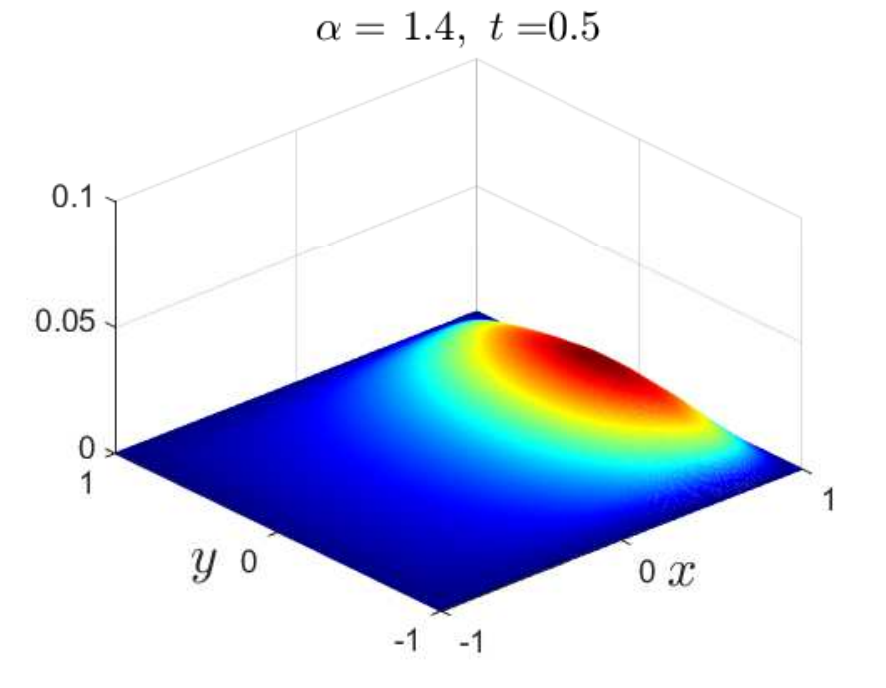}\hspace{-5mm}
\includegraphics[height = 3.96cm, width = 4.36cm]{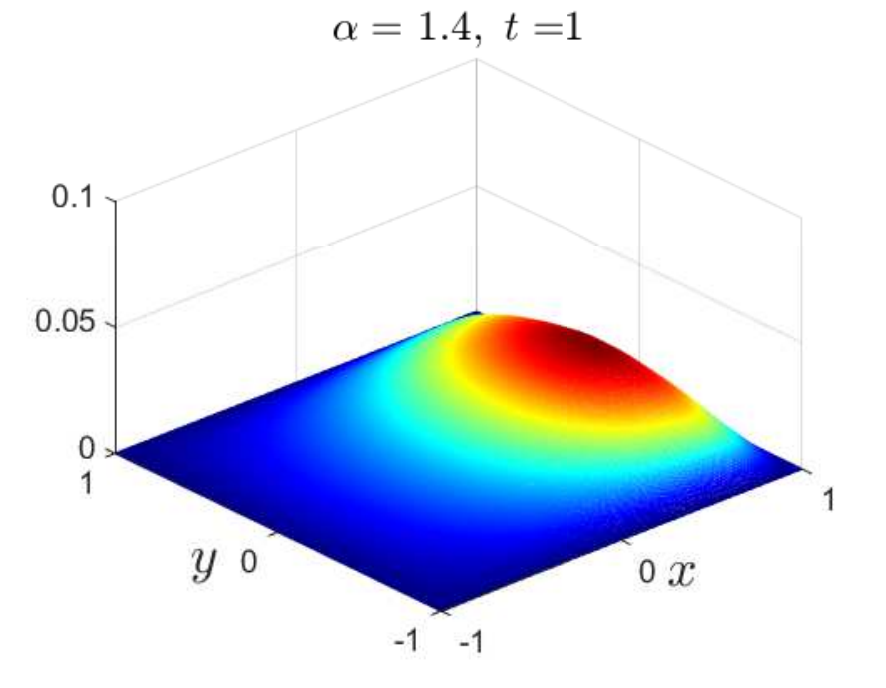}\hspace{-5mm}
\includegraphics[height = 3.96cm, width = 4.36cm]{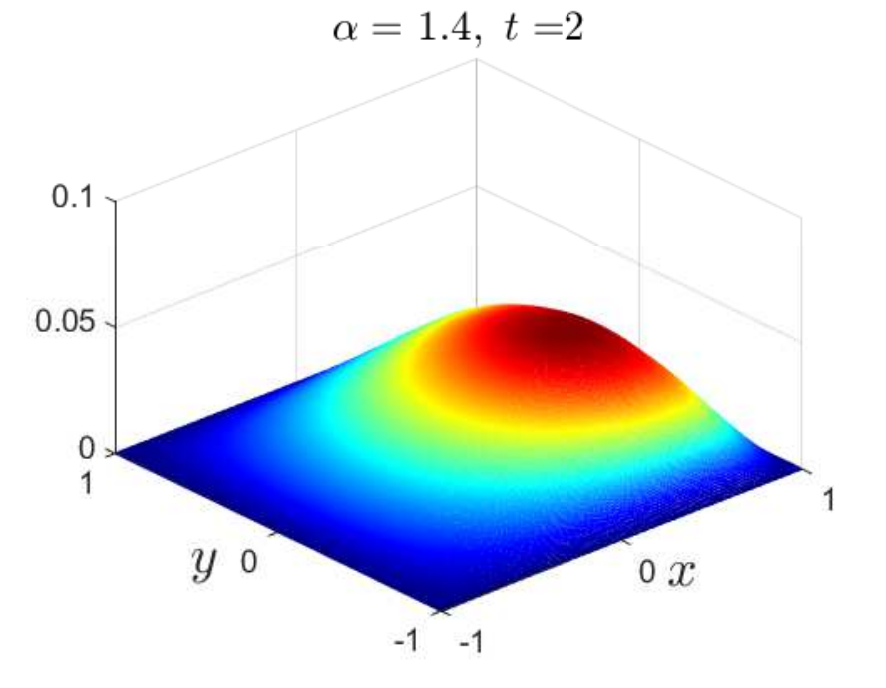}\hspace{-5mm}
\includegraphics[height = 3.96cm, width = 4.36cm]{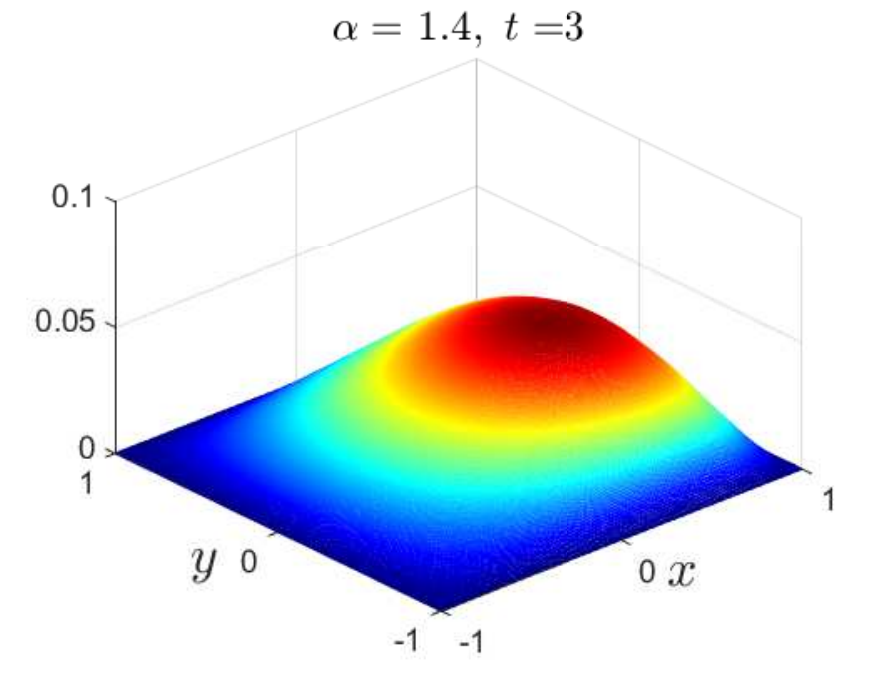}}
\caption{Numerical solution of heat equation \eqref{fun4-2-1} at different time, where $x_c = 1.3$. } \label{Fig5-3-3}
\end{figure}

Next, we continue our study of boundary effects in Fig. \ref{Fig5-3-3},  for $x_c = 1.3$ and $\ap = 0.7,\, 1.4$. 
In this case, the nonzero boundary conditions are imposed on region $\Xi = [1.3, 1.75]\times [-1, 1]$ which has no contact with the domain $\Og$. 
Hence, the solution of classical ($\ap = 2$) problems remain $u(\bx, t) \equiv 0$ for time $t \ge 0$, due to the homogeneous boundary conditions on $\p\Og$. 
In contrast to classical cases, the nonzero boundary conditions on region $\Xi$  lead to a growth of solution in the fractional cases, although there is no contact between $\Xi$ and $\Og$.  
This observation suggests the main difference between classical and fractional Laplacians -- one is local, and the other is nonlocal.  
The fractional Laplacian describing the L\'evy flights can transfer the boundary information into the domain, even though $\Xi$ and $\Og$ are not connected.

Fig. \ref{Fig5-3-3} further shows that the smaller the power $\ap$, the stronger the boundary effects. 
Comparing the results of $\ap = 1.4$ in Figs. \ref{Fig5-3-1} and \ref{Fig5-3-3}, we find  that the boundary effects are stronger if the distance between  $\Xi$ and $\Og$ (i.e., the distance is $x_c - 1$) is smaller. 
These two phenomena can be explained by the kernel function $|\bx - \by|^{-(d+\ap)}$ of the fractional Laplacian in \eqref{integralFL}, whose value is affected by both distance $|\bx - \by|$ and power $\ap$. 
The smaller the power $\ap$ (or the shorter the distance $|\bx-\by|$), the larger the kernel function, and thus the stronger the interactions from boundary conditions. 
Our study might provide insights for boundary control in the study of fractional PDEs \cite{Antil2019, Delia2019}. 

\section{Conclusion}
\label{section6}

We proposed a new meshless pseudospectral method based on the generalized inverse multiquadric (GIMQ) functions to solve problems with $(-\Dt)^{\fl{\ap}{2}}$ for $0 < \ap \le 2$. 
The operator $(-\Dt)^{\fl{\ap}{2}}$ represents the local classical Laplacian $-\Dt$ when $\ap = 2$,  while  the nonlocal fractional Laplacian if $\ap < 2$. 
In existing literature, numerical discretizations of classical and fractional Laplacians were carried out separately, due to their different nature. 
Consequently, their computer implementations were usually incompatible. 
On the other hand, it is well-known that the fractional Laplacian $(-\Dt)^\fl{\ap}{2}$ can be defined via the pseudodifferential operator in \eqref{pseudo} or the pointwise hypersingular integral form in \eqref{integralFL}.  
We combined the advantages of these two definitions and proposed an $\ap$-parametric method for the Laplacian $(-\Dt)^{\fl{\ap}{2}}$ with $\ap \in (0, 2]$. 
Hence, our method unifies the discretization of classical and fractional Laplacians. 
Moreover, it can easily incorporate Dirichlet boundary conditions.

In our method, we used the GIMQ functions  with power $\bt = -(d+1)/2$ (also known as the Poisson kernel), for dimension $d \ge 1$. 
Notice that the Laplacian of GIMQ functions can be analytically written by the hypergeometric function for any $\ap \ge 0$. 
This provides the key to avoiding numerical evaluations to the hypersingular integral and unifying approximations of classical and fractional Laplacians. 
Our studies of the diffusion problems illustrated the differences between the normal  ($\ap = 2$) and anomalous ($\ap < 2$) diffusion. 
It showed that boundary conditions of anomalous diffusion problems can affect the solution in domain $\Og$ even though they may not have direct contact, in contrast to the normal diffusion where the boundary conditions have to be imposed on $\p\Og$. 
These studies could provide insights of simulations and applications with the fractional models.

We compared our method with the recent Gaussian-based method in \cite{Burkardt0020} and found that both methods have the spectral accuracy.  
Furthermore, we found that if uniformly distributed RBF center points are considered,  the optimal shape parameter of Gaussian-based method is more sensitive to the number of points. 
This suggests that different strategies might be considered in selecting the shape parameters of Gaussian and GIMQ based methods. 
Two approaches in selecting shape parameters of our GIMQ-based method were studied in detail.  
Numerical studies showed that the condition-number indicated shape parameter can effectively suppress the ill-conditioning issues and thus guarantee the numerical accuracy when a larger number of points are used. 
This method may take extra computational time in searching a proper shape parameter to satisfy the desired condition number range. 
While the random-perturbed shape parameters can save the time of searching. \\

\bigskip
\noindent{\bf Acknowledgements. }
This work was supported by the US National Science Foundation under Grant Number DMS-1913293 and DMS-1953177.

\end{document}